\documentclass[a4paper,12pt]{amsart}
\usepackage[utf8]{inputenc}

\usepackage{amsmath,amssymb,amsthm}
\usepackage{bbm}
\usepackage{bm}
\usepackage{xspace}
\usepackage{csquotes}
\usepackage{standalone}

\usepackage[x11names,table]{xcolor}
\usepackage[naturalnames=false]{hyperref}
\usepackage{graphicx}
\graphicspath{ {img/} }
\usepackage{tikz}
\usetikzlibrary{calc}
\usetikzlibrary{cd}
\usetikzlibrary{shapes}
\usepackage{csquotes}
\usepackage{booktabs}
\hypersetup{
	colorlinks = true,          
	urlcolor = DeepSkyBlue3, linkcolor = Chartreuse4, citecolor = DarkOrange2
}
\usepackage{fullpage}
\RequirePackage[backend=bibtex,isbn=false,url=false,bibencoding=utf8,maxbibnames=4]{biblatex}
\renewbibmacro{in:}{, }
\DeclareFieldFormat[misc]{title}{\mkbibquote{#1}}
\DeclareFieldFormat[article]{volume}{\mkbibbold{#1}}
\DeclareFieldFormat{url}{\textsc{url}: \href{#1}{#1}}
\DeclareFieldFormat{doi}{\textsc{doi}: \href{http://dx.doi.org/#1}{#1}}
\DeclareFieldFormat{eprint:arxiv}{arXiv:\href{https://arxiv.org/abs/#1}{#1}}
\usepackage{pgfplots}
\pgfplotsset{compat=1.5}

\usepackage{mathbbol}

\newcommand\PP{{\mathbb P}}

\newcommand\RR{{\mathbb R}}

\newcommand\TT{{\mathbb T}}

\newcommand\cE{{\mathcal E}}
\newcommand\cS{{\mathcal S}}

\newcommand\cH{{\mathcal H}}

\newcommand\cT{{\mathcal T}}

\DeclareMathSymbol{\0}{\mathord}{bbold}{`0}
\DeclareMathSymbol{\1}{\mathord}{bbold}{`1}

\newcommand\SetOf[2]{\left\{#1 \mid #2\right\}}

\newcommand\bigSetOf[2]{\bigl\{#1 \bigm| #2\bigr\}}

\newcommand\torus[1]{\RR^{#1}/\RR\1} 

\newcommand\TTmax{\TT^{\max}}
\newcommand\TTmin{\TT^{\min}}
\newcommand\TP{{\TT\PP}}
\newcommand\BergmanFan{{\widetilde B}}

\DeclareMathOperator\dist{dist}
\DeclareMathOperator\conv{conv}
\DeclareMathOperator\tconv{tconv}

\DeclareMathOperator\codim{codim}
\DeclareMathOperator\lcm{lcm}

\DeclareMathOperator{\FW}{FW} 
\newcommand{\FWsym}{\FW_{\rm sym}} 
\DeclareMathOperator\CovDec{CovDec} 
\DeclareMathOperator\TGr{TGr} 

\newcommand\transpose[1]{#1^\top}

\newtheorem{theorem}{Theorem}

\newtheorem{proposition}[theorem]{Proposition}
\newtheorem{corollary}[theorem]{Corollary}

\theoremstyle{remark}

\newtheorem{example}[theorem]{Example}
\newtheorem{definition}[theorem]{Definition}
\newtheorem{remark}[theorem]{Remark}

\newcommand\polymake{\texttt{polymake}\xspace}

\newcommand\mcf{\texttt{mcf}\xspace}

\definecolor{tolorange}{RGB}{238,119,51}  
\definecolor{tolblue}{RGB}{0,119,187}     
\definecolor{tolgreen}{RGB}{0,153,136}    
\definecolor{tolpurple}{RGB}{136,34,85}   
\definecolor{ibmyellow}{RGB}{255, 176, 0} 
\definecolor{darkgray}{RGB}{64,64,64} 

\tikzstyle{blackdot} = [circle,draw=black,fill=black,scale=0.5]

\usepackage{expl3}
\ExplSyntaxOn
\newcommand \breakDOI[1]
 {
  \tl_map_inline:nn { #1 } { \href{http://dx.doi.org/#1}{##1} \penalty0 \scan_stop: }
 }
\ExplSyntaxOff

\DeclareFieldFormat{doi}{%
  \mkbibacro{DOI}\addcolon\space
  {\breakDOI{#1}}}

\title{Tropical medians by transportation}

\author{Andrei Com\u{a}neci \and Michael Joswig}

\address[Andrei Com\u{a}neci]{
	Technische Universität Berlin,
	Chair of Discrete Mathematics/Geometry \\
	\texttt{comaneci@math.tu-berlin.de}	
}

\address[Michael Joswig]{
	Technische Universität Berlin,
	Chair of Discrete Mathematics/Geometry \\
	Max-Planck Institute for Mathematics in the Sciences, Leipzig \\
	\texttt{joswig@math.tu-berlin.de}
}

\thanks{Support by the Deutsche Forschungsgemeinschaft (DFG, German Research Foundation) under Germany's Excellence Strategy -- The Berlin Mathematics Research Center MATH$^+$ (EXC-2046/1, project ID 390685689) gratefully acknowledged.
  M.~Joswig has further been supported by \enquote{Symbolic Tools in Mathematics and their Application} (TRR 195, project-ID 286237555).
  A.~Com\u{a}neci was supported by \enquote{Facets of Complexity} (GRK 2434, project-ID 385256563).
		}

\begin{document}

\begin{abstract}
  Fermat--Weber points with respect to an asymmetric tropical distance function are studied.
  It turns out that they correspond to the optimal solutions of a transportation problem.
  The results are applied to obtain a new method for computing consensus trees in phylogenetics.
  This method has several desirable properties; e.g., it is Pareto and co-Pareto on rooted triplets.
\end{abstract}

\maketitle

\section{Introduction} \label{sec:intro}
\noindent
The following optimization problem can be studied in any metric space.
Given a finite number of points, sometimes called \emph{sites}, find a point which minimizes the sum of the distances to the sites.
Such a point is called a \emph{Fermat--Weber} point, and this is some version of a geometric median of the sites, which is known to be robust in a certain sense \cite[\textsection 21]{BMS:1999}. 
Computing Fermat--Weber points is a rich topic with a remarkable history; see \cite[Chapter~II]{BMS:1999}. 
Here we consider a specific distance function, $d_\triangle$, which occurs in tropical geometry \cite{RRlattice,LapLatticeGraph}. 
This function is asymmetric, i.e., $d_\triangle(a,b)$ may differ from $d_\triangle(b,a)$.
So we call $d_\triangle$ the \emph{asymmetric tropical distance}, to differentiate from the \emph{symmetric tropical distance}, which is more common \cite[\textsection 5.3]{ETC}.
The symmetric tropical Fermat--Weber problem was investigated by Lin, Sturmfels, Tang and Yoshida \cite{ConvTreeSp} and Lin and Yoshida~\cite{Lin-Yoshida:2018}.

As our main result we prove that one (asymmetric tropical) Fermat--Weber point can be computed by solving a transportation problem.
The transportation problem is an optimization classic, with numerous applications, both in theory and practice.
For an overview we refer to Schrijver's monograph \cite[\textsection 21.6]{Schrijver03:CO_A}; see also the survey by De Loera and Kim for the polyhedral geometry point of view \cite{DeLoera+Kim:2014}.
Efficient methods for solving transportation problems include algorithms by Tokuyama and Nakano \cite{Tokuyama+Nakano:1995}, Kleinschmidt and Schannath \cite{Kleinschmidt+Schannath:1995}, and Brenner \cite{Brenner:2008}.
In general, Fermat--Weber points are not unique, so one part of the present work is devoted to understanding the entire (asymmetric tropical) \emph{Fermat--Weber set} for a given set of sites.
This is tightly related to the study of tropical hyperplane arrangements and tropical convex hulls, which are at the core of tropical combinatorics \cite{ETC}.
The latter subject is concerned with the rich interplay between tropical geometry and optimization.

One motivation for studying the Fermat--Weber problem in the setting of tropical geometry comes from phylogenetics \cite{ConvTreeSp,Lin-Yoshida:2018}.
In that field, a part of computational biology, the goal is to associate trees to input data.
A typical example are trees encoding ancestral relations among species, where the data originates from strands of DNA.
In tropical geometry spaces of metric trees with $n$ labeled leaves occur naturally as the tropical Grassmannians $\TGr(2,n)$; see \cite[\textsection 4.3]{Tropical+Book} and \cite[\textsection 10.6]{ETC}.
In phylogenetics many different methods are known to construct trees from a fixed data set.
Since those methods usually do not come up with the same tree, there is a need to find the common ground.
This gives rise to some \emph{consensus tree}, which describes where the several methods agree.
Finding consensus trees is a topic of its own \cite{Bryant:2003}, and the authors of \cite{ConvTreeSp} argue that \enquote{tropical convexity and tropical linear algebra $\ldots$ behave better} than other methods.
Here we show that passing from the symmetric tropical distance function to its asymmetric sibling leads to a new method for computing metric consensus trees which is even better behaved.
This is because the asymmetric tropical Fermat--Weber sets are nicer geometrically.
In particular, we show that our approach leads to a consensus method which is regular in the sense of \cite{Bryant+Francis+Steel:2017}.
Such a procedure is not known for symmetric tropical distances.
For the purpose of finding a consensus method in tree space, there is no immediate disadvantage of employing an asymmetric distance rather than a symmetric one.
In fact, asymmetric distances are common in location theory \cite{NickelPuerto:2005}.

Our paper is organized as follows.
We start out with a brief summary of facts from tropical combinatorics which are relevant for the Fermat--Weber problem.
Then we prove that the Fermat--Weber set arises as a cell in the covector decomposition of the tropical torus induced by the sites.
That covector cell is then characterized in several ways.
The first approach employs regular subdivisions of products of simplices; see \cite[\textsection 6.2]{DLRS}.
That leads to a linear programming formulation of the Fermat--Weber problem, and the dual linear program is a transportation problem.
The latter then provides efficient algorithms.
The final chapter is devoted to computing Fermat--Weber sets in spaces of equidistant trees.
In addition to theoretical results we report on computational experiments with \polymake \cite{DMV:polymake} and \mcf \cite{mcf}.

\noindent \textbf{Related work.}
As an early contribution of tropical geometry to data science G\"{a}rtner and Jaggi \cite{GaertnerJaggi:2006} developed a concept for \enquote{tropical support vector machines}, with applications to classification in mind.
A different train of thought was developed by Pachter and Sturmfels \cite[\textsection 2.4]{ASCB} who related tropical geometry to phylogenetic trees.
Later, Lin et al. \cite{ConvTreeSp} connected these ideas to the geometry of tree spaces studied by Billera, Holmes and Vogtman \cite{BilleraHolmesVogtmann01}.
Yoshida, Zhang and Zhang \cite{YoshidaZhangZhang:2019} proposed a method to analyze data, which they call \enquote{tropical principal component analysis}.
In a way, the latter may be viewed as a synthesis of the above.
A key contribution here is work of Ardila and Klivans, who saw that spaces of equidistant trees form the Bergman fans of the graphic matroids of complete graphs \cite{ArdilaKlivans:2006}.
The term \enquote{tropical convexity} was coined by Develin and Sturmfels \cite{Develin+Sturmfels:2004} to connect tropical geometry with the older topic of $(\max,+)$-linear algebra \cite{Butkovic:2010}.

\section{Tropical convexity} \label{sec:tconv}
\noindent
The purpose of this section is to set our notation and to collect a few facts which are key to our methods; for the details we refer to \cite{ETC}.
We consider the \emph{tropical semiring} $\TTmax=(\RR\cup\{-\infty\},\oplus,\odot)$ with $\oplus=\max$ as the tropical addition and $\odot=+$ as the tropical multiplication.
The additive neutral element is $-\infty$, and $0$ is neutral with respect to the multiplication.
Usually, we abbreviate $\TT=\TTmax$.
The set $\TT^n$ inherits the structure of a semimodule by componentwise tropical addition and tropical scalar multiplication.

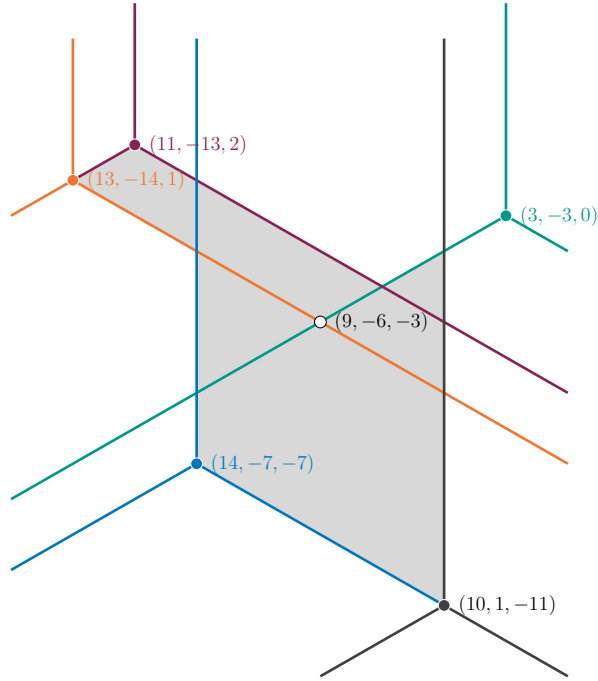
\begin{figure}[th]
  \centering
	\resizebox{8cm}{!}{
	\begin{tikzpicture}[scale=1.5,ultra thick
,x={(0.866,-0.5)},y={(0,1)},z={(-0.866,-0.5)}]
		
		
                \fill[gray!60,fill opacity=0.5] (-9,-4) -- (-7,-4) -- (-7,-7) -- (-3,-7) -- (-3,-2) -- (-2,-1) -- (-4,-3) -- (-8,-3) -- cycle;
                
		\draw[color = tolgreen] (-10,-9) -- (-2,-1) -- (-2,2);
		\draw[color = tolgreen] (-2,-1) -- (-1,-1);
		\draw[color = tolorange] (-10,-5) -- (-9,-4) -- (-9,-2);
		\draw[color = tolorange] (-9,-4) -- (-1,-4);
		\draw[color = tolpurple] (-9,-4) -- (-8,-3) -- (-8,-1);
		\draw[color = tolpurple] (-8,-3) -- (-1,-3);
		\draw[color = darkgray] (-5,-9) -- (-3,-7) -- (-3,1);
		\draw[color = darkgray] (-3,-7) -- (-1,-7);
		\draw[color = tolblue] (-10,-10) -- (-7,-7) -- (-7,-1);
		\draw[color = tolblue] (-7,-7) -- (-3,-7);
		
		
		\node at (-2,-1) [circle, fill = tolgreen, draw = white, semithick, inner sep = 2.5pt, label={[tolgreen]right:$(3,-3,0)$}] {};
		\node at (-9,-4) [circle, fill = tolorange, draw = white, semithick, inner sep = 2.5pt, label={[tolorange]right:$(13,-14,1)$}] {};
		\node at (-8,-3) [circle, fill = tolpurple, draw = white, semithick, inner sep = 2.5pt, label={[tolpurple]right:$(11,-13,2)$}] {};
		\node at (-3,-7) [circle, fill = darkgray, draw = white, semithick, inner sep = 2.5pt, label={[black]right:$(10,1,-11)$}] {};
		\node at (-7,-7) [circle, fill = tolblue, draw = white, semithick, inner sep = 2.5pt, label={[tolblue]right:$(14,-7,-7)$}] {};
		
		\node (FW) at (-5,-4) [circle, draw= black, fill = white, semithick, inner sep = 2.5pt, label={[black]right:$(9,-6,-3)$}] {};
		
\end{tikzpicture}
	}
  \caption{Five sites in $\torus{3}$, their $\max$-tropical convex hull, the induced $\min$-tropical hyperplane arrangement, and the unique Fermat--Weber point (marked white)}
  \label{fig:FWpoint}
\end{figure}

A \emph{tropical cone} in $\TT^n$ is a nonempty subset $C$ which contains each tropical linear combination $\lambda {\odot} x \oplus \mu {\odot} y$ for $\lambda,\mu\in\TT$ and $x,y\in C$.
Each tropical cone contains the point $-\infty\1$ and the entire set $\RR\1$, where $\1$ is the all-ones vector.
Therefore it is convenient to pass to the quotient $\TP^{n-1}:=(\TT^n\setminus\{-\infty\1\})/\RR\1$, which is called the \emph{tropical projective space}.
A subset of $\TP^{n-1}$ is \emph{tropically convex} if it arises as the image of a tropical cone under the canonical projection.
A \emph{tropical polytope} is a finitely generated tropically convex set.
The \emph{tropical projective torus} $\torus{n}$ is the subset of points in $\TP^{n-1}$ with finite coordinates.
We say that a tropical polytope is \emph{bounded} if it lies in $\torus{n}$.

Tropical convexity is intimately related to ordinary convexity, polyhedral geometry and (linear) optimization.
For instance, tropical polytopes arise as the images of ordinary convex polytopes over ordered fields of real Puiseux series under the valuation map; see \cite[Observation~5.10]{ETC}.
Yet, here the following less algebraic description is more relevant.

We consider an arbitrary $m{\times}n$-matrix $V=(v_{ij})$ with real coefficients.
Each row $v_i=(v_{i1},\dots,v_{in})$ is a point in $\RR^n$ (or $\torus{n}$, if we ignore the tropical scaling).
So $V$ may be viewed as a configuration of $m$ labeled points in $\torus{n}$.
Technically, it is convenient to assume that each ordinary row sum equals zero, i.e., each row lies in the set
\[
  \cH \ = \ \SetOf{x\in\RR^n}{x_1+x_2+\dots+x_n=0} \enspace.
\]
Observe that each point in $\torus{n}$ has a unique representative in $\RR^n$ which lies in $\cH$.
So we can identify the tropical projective torus $\torus{n}$ with the ordinary linear hyperplane $\cH$ in $\RR^n$.
This also works topologically, since the quotient vector space $\torus{n}$ is homeomorphic to $\RR^{n-1}$ (and thus with $\cH$ considered as a subset of $\RR^n$).
The tropical projective space $\TP^{n-1}$ is a compactification of $\torus{n}$, where the boundary comprises those points which have at least one infinite coordinate.

Adding vectors tropically works coefficient-wise, and there is also tropical multiplication by scalars.
With these notions, the \emph{$\max$-tropical convex hull} of (the rows of) $V$ is
\[
  \tconv^{\max}(V) \ := \ \bigSetOf{ \lambda_1{\odot} v_1 \oplus \dots \oplus \lambda_m{\odot}v_m}{\lambda_i\in\RR} + \RR\1 \enspace,
\]
which is a subset of $\torus{n}$.
The rows of the matrix $V$ also define an arrangement of $m$ tropical hyperplanes in $\torus{n}$, with respect to $\min$ as the tropical addition, and we denote this $\cT(V)$.
In this context each row arises as the \emph{apex} of a $\min$-tropical hyperplane.
Everything that we explained above also works for the $\min$-tropical semiring $\TTmin=(\RR\cup\{\infty\},\min,+)$, which is isomorphic to $\TTmax$ as a semiring via $x\mapsto -x$.
Observe that this involution leaves the hyperplane $\cH$ invariant.

The $\min$-tropical hyperplane arrangement with the rows of $V$ as their apices induces an ordinary polyhedral subdivision, $\CovDec(V)$, of $\torus{n}$ (or, equivalently, $\cH$), called the \emph{covector subdivision} induced by (the rows of) $V$.
Its cells, which are called \emph{covector cells}, are convex in three different senses: they are $\max$-tropically convex, $\min$-tropically convex and convex in the ordinary sense (as subsets of $\cH$).
Such polyhedra are known as \emph{polytropes}, and they may be bounded or unbounded.
As ordinary polyhedra, the polytropes are characterized by the property that their facet normal directions are $e_i-e_j$ for $i,j\in[n]$ distinct.
The union of the bounded covector cells equals the tropical convex hull $\tconv^{\max}(V)$.

\begin{figure}[th]
  \centering
	\resizebox{8cm}{!}{\tikzset{Hypersurface/.style = {ultra thick}}
\tikzset{InnerTriangle/.style = {fill=gray!50, draw=gray, thick}}
\tikzset{Triangle/.style = {ultra thick, black}}

\begin{tikzpicture}[scale = 1.0]

\coordinate (B) at (-210:4cm);
\coordinate (A) at (30:4cm);
\coordinate (C) at (-90:4cm);

\newcommand\bc[3]{(barycentric cs:A=#1,B=#2,C=#3)}

\filldraw[InnerTriangle] \bc{5}{0}{0} -- \bc{4}{1}{0} -- \bc{4}{0}{1} -- cycle;
\filldraw[InnerTriangle] \bc{2}{3}{0} -- \bc{1}{4}{0} -- \bc{1}{3}{1} -- cycle;
\filldraw[InnerTriangle] \bc{1}{4}{0} -- \bc{1}{3}{1} -- \bc{0}{3}{2} -- \bc{0}{5}{0} -- cycle;
\filldraw[InnerTriangle] \bc{0}{2}{3} -- \bc{1}{1}{3} -- \bc{0}{1}{4} -- cycle;
\filldraw[InnerTriangle] \bc{0}{1}{4} -- \bc{1}{1}{3} -- \bc{2}{0}{3} -- \bc{0}{0}{5} -- cycle;

\draw[InnerTriangle,fill=none] \bc{0}{3}{2} -- \bc{1}{2}{2} -- \bc{1}{1}{3} -- \bc{2}{1}{2} -- \bc{3}{0}{2};
\draw[InnerTriangle,fill=none] \bc{3}{2}{0} -- \bc{3}{1}{1} -- \bc{2}{1}{2} -- \bc{2}{2}{1} -- \bc{1}{2}{2};
\draw[InnerTriangle,fill=none] \bc{1}{3}{1} -- \bc{2}{2}{1} -- \bc{3}{2}{0};
\draw[InnerTriangle,fill=none] \bc{3}{1}{1} -- \bc{4}{0}{1};

\draw[tolpurple,Hypersurface] \bc{4}{10}{1} -- \bc{3}{6}{1} -- \bc{5}{4}{1} -- \bc{5}{2}{3} -- \bc{7}{0}{3};
\draw[tolpurple,Hypersurface] \bc{4}{10}{1} -- \bc{3}{7}{0};
\draw[tolpurple,Hypersurface] \bc{4}{10}{1} -- \bc{2}{15}{3};
\draw[tolblue,Hypersurface] \bc{1}{4}{10} -- \bc{0}{3}{7};
\draw[tolblue,Hypersurface] \bc{1}{4}{10} -- \bc{1}{3}{6} -- \bc{1}{5}{4} -- \bc{3}{5}{2} -- \bc{5}{5}{0};
\draw[tolblue,Hypersurface] \bc{1}{4}{10} -- \bc{3}{2}{15};
\draw[tolgreen,Hypersurface] \bc{13}{1}{1} -- \bc{9}{1}{0};
\draw[tolgreen,Hypersurface] \bc{13}{1}{1} -- \bc{9}{0}{1};
\draw[tolgreen,Hypersurface] \bc{13}{1}{1} -- \bc{8}{1}{1} -- \bc{6}{3}{1} -- \bc{4}{3}{3} -- \bc{2}{3}{5} -- \bc{0}{5}{5};
\draw[darkgray,Hypersurface] \bc{3}{2}{15} -- \bc{1}{0}{4};
\draw[darkgray,Hypersurface] \bc{3}{2}{15} -- \bc{0}{1}{9};
\draw[darkgray,Hypersurface] \bc{3}{2}{15} -- \bc{3}{1}{6} -- \bc{5}{1}{4} -- \bc{7}{1}{2} -- \bc{7}{3}{0};
\draw[tolorange,Hypersurface] \bc{2}{15}{3} -- \bc{0}{4}{1};
\draw[tolorange,Hypersurface] \bc{2}{15}{3} -- \bc{1}{9}{0};
\draw[tolorange,Hypersurface] \bc{2}{15}{3} -- \bc{1}{6}{3} -- \bc{3}{4}{3} -- \bc{3}{2}{5} -- \bc{5}{0}{5};

\draw[Triangle] (A) -- (B) -- (C) -- cycle;

\foreach \i in {0,...,5} {
   \pgfmathsetmacro{\ci}{5-\i}
   \foreach \j in {0,...,\ci} {
     \pgfmathsetmacro{\k}{int(5-\i-\j)}
     \node[fill=white,shape=ellipse,inner sep=1pt] (\i\j) at \bc{\i}{\j}{\k} {$\i \j
\k$};
   }
}

\end{tikzpicture}	}
  \caption{%
    Mixed subdivision $\cS(V)$ of $5\cdot\Delta_2$ for $V$ as in Example~\ref{exmp:FWpoint}. 
    The 21 lattice points of $5\cdot\Delta_2$ are marked with their coordinates.
    The $\min$-tropical hyperplane arrangement $\cT(V)$ admits a piecewise-linear embedding}
  \label{fig:dualSubd}
\end{figure}
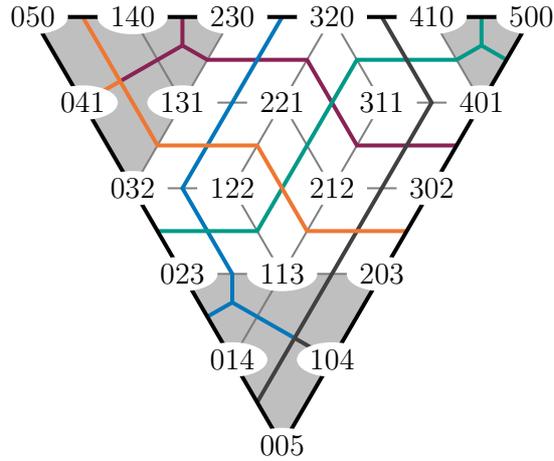

\begin{example}\label{exmp:FWpoint}
  We illustrate the various concepts from tropical convexity for the matrix $V\in\RR^{5\times3}$ whose transpose reads
  \[
    \transpose{V} \ = \ \begin{pmatrix}
      14  & 13  & 11  & 10  & 3 \\
      -7  & -14 & -13 & 1   & -3  \\
      -7 & 1   & 2   & -11 & 0  \\
    \end{pmatrix} \enspace .
  \]
  The rows of $V$ (equivalently, the columns of $\transpose{V}$) encode five points in $\torus{3}$; see Fig.~\ref{fig:FWpoint}.
  The covector decomposition has 15 regions of maximal dimension $2$, and six of them are bounded.
  The $\max$-tropical convex hull $\tconv^{\max}(V)$ consists of the union of bounded cells, which are shaded gray; it also contains the green line segment extending from $(3,-3,0)$ to the lower left.
\end{example}

We consider the \emph{envelope}
\begin{equation}\label{eq:envelope}
  \cE(V) \ := \ \bigSetOf{(t,x)\in\RR^m\times\RR^n}{t_i+x_j \geq v_{ij}} \enspace ,
\end{equation}
which is an unbounded ordinary polyhedron.
By \cite[Theorem 6.14]{ETC} the cells of $\CovDec(V)$ arise as the images of faces of $\cE(V)$ under the coordinate projection $(t,x) \mapsto x$.
Moreover, $\CovDec(V)$ is dual to the regular subdivision, $\Sigma(V)$, of the product of simplices $\Delta_{m-1}\times\Delta_{n-1}=\conv\SetOf{(e_i,e_j)}{i\in[m],\, j\in[n]}$ obtained from lifting $(e_i,e_j)$ to the height $v_{ij}$.
Here we take regular subdivisions induced by upper convex hulls since $\max$ is our tropical addition.
For the same reason the inequality sign \enquote{$\geq$} is reversed in comparison with the $\min$-tropical version in \cite[(6.1)]{ETC}.
Via the Cayley trick the subdivision $\Sigma(V)$ of $\Delta_{m-1}\times\Delta_{n-1}$ corresponds to a mixed subdivision, $\cS(V)$, of the dilated simplex $m\cdot\Delta_{n-1}$; see \cite[\textsection 9.2]{DLRS} and \cite[\textsection 4.5]{ETC}.
For instance, this is convenient for properly visualizing $\Sigma(V)$, which is a polyhedral complex of dimension $(m-1)(n-1)$.
Fig.~\ref{fig:dualSubd} shows the mixed subdivision $\cS(V)$ of $5\cdot\Delta_2$ for the matrix $V$ from Example~\ref{exmp:FWpoint}.

\section{Fermat--Weber sets} \label{sec:fermat-weber}
\noindent
We examine the Fermat--Weber problem through tropical combinatorics and polyhedral geometry.
As our first key observation we show that asymmetric tropical Fermat--Weber sets arise as specific cells in the covector subdivisions induced by the sites.

The \emph{asymmetric tropical distance} in $\RR^n$ is given by
\begin{equation}\label{eq:d-triangle}
  d_\triangle(a,b) \ = \ \sum_{i\in[n]} (b_i-a_i) - n \min_{i\in[n]}(b_i-a_i) \ = \ \sum_{i\in[n]} (b_i-a_i) + n \max_{i\in[n]}(a_i-b_i) \enspace ,
\end{equation}
where $a, b\in\RR^n$.
Since $d_\triangle(a',b')=d_\triangle(a,b)$ for $a-a'\in\RR\1$ and $b-b'\in\RR\1$, this induces a directed distance function in the $(n{-}1)$-dimensional quotient vector space $\torus{n}$.
We do not distinguish between the function $d_\triangle:\RR^n\times\RR^n\to\RR_ {\geq0}$ and the induced function on $\torus{n}$.
The asymmetric tropical distance is a \enquote{polyhedral norm} with respect to the standard simplex $\Delta_{n-1}:=\conv\{e_1,\dots,e_n\}+\RR\1$ in $\torus{n}$; see \cite[\textsection 20]{BMS:1999}. 
This may also be seen as a rescaled version of the \enquote{tropical Funk metric} discussed in \cite[\textsection 3.3]{ABGJ:2018}.
More common in tropical geometry is the \emph{symmetric tropical distance} between $a,b\in\RR^n$ (or $\torus{n}$).
It is defined as
\[
  \dist(a,b) \ = \ \max_{i\in[n]} \left(a_i-b_i\right) - \min_{j\in[n]} \left(a_j-b_j\right) \ = \ \max_{i,j\in[n]}(a_i-b_i-a_j+b_j) \enspace ;
\]
see \cite[\textsection 5.3]{ETC}.
We have $\dist(a,b) = \tfrac{1}{n}(d_{\triangle}(a,b) + d_\triangle(b,a))$.
Throughout this section we fix a finite set of sites $V=\{v_1,v_2,\dots,v_m\}$ in $\cH$, which we identify with $\torus{n}$.
\begin{definition}
  An \emph{(asymmetric tropical) Fermat--Weber point} with respect to $V$ is a point in $\cH$ which minimizes the sum of the asymmetric tropical distances from these sites.
\end{definition}
In general such a point is not unique.
Hence, we let $\FW(V)$ denote the set of all asymmetric tropical Fermat--Weber points and call it the \emph{(asymmetric tropical) Fermat--Weber set} with respect to $V$.
This is the asymmetric analog of the symmetric tropical Fermat--Weber set studied in \cite{Lin-Yoshida:2018}.

Fixing the site $v_i\in V$, the distance function from $v_i$, which reads
\begin{equation}\label{eq:from-the-sites}
  d_\triangle(v_i,x) \ = \ n \max_{j\in[n]} (v_{ij}-x_j) \quad \text{for } x\in\cH \enspace,
\end{equation}
is convex in the ordinary sense and piecewise linear.
Its regions of linearity are precisely the $n$ closed sectors of the $\min$-tropical hyperplane with apex $v_i$; see \cite[\textsection 5.5]{ETC}.
Consequently, the common subdivision of the regions of linearity of the sites is exactly the covector decomposition $\CovDec(V)$.
Our first main theorem shows that the Fermat--Weber set $\FW(V)$ is a bounded cell of that subdivision.

For the sake of a precise formulation of that result, we pass to the regular triangulation, $\Sigma(V)$, of $\Delta_{m-1}\times\Delta_{n-1}$ which is dual to the covector subdivision $\CovDec(V)$.
The relatively open cells of $\Sigma(V)$ partition the product of the ordinary polytope $\Delta_{m-1}\times\Delta_{n-1}$.
The point $(\tfrac{1}{m}\1,\tfrac{1}{n}\1)$ is the vertex barycenter of  $\Delta_{m-1}\times\Delta_{n-1}$.
So there is a unique cell, $C_V$, which contains that $(\tfrac{1}{m}\1,\tfrac{1}{n}\1)$ in its relatively interior.
This will play an important role in our study of Fermat--Weber points.

\begin{definition}
We call the unique cell of $\Sigma(V)$ which contains that $(\tfrac{1}{m}\1,\tfrac{1}{n}\1)$ in its relatively interior the \emph{central cell} of $\Sigma(V)$, and denote it by $C_V$.
Its dual in $\CovDec(V)$ will be called the \emph{central covector cell}.
\end{definition}

\begin{theorem} \label{thm:convFW}
  The Fermat--Weber set $\FW(V)$ agrees with the central covector cell in $\CovDec(V)$.
  In particular, $\FW(V)$ is a bounded polytrope in $\cH$, and it is contained in the tropical polytope $\tconv^{\max}(V)$.
\end{theorem}
\begin{proof}
  Consider the linear program
  \begin{equation} \label{lp:prim}
    \begin{array}{ll@{}ll}
      \text{minimize}  & \displaystyle  n \cdot (t_{1} + \dots +  t_{m}) & \\
      \text{subject to}& \displaystyle v_{ij} - x_j \ \leq \ t_i \,, & \quad \text{for } i\in[m] \text{ and } j\in[n] \\
                       & \displaystyle x_{1} + \dots + x_{n} \ = \ 0  &
    \end{array}
  \end{equation}
  with $mn+1$ constraints in the $m+n$ variables $t_1,t_2,\dots,t_m,x_1,x_2,\dots, x_n$.
  The coefficients $v_{ij}$ are the coordinates of the sites.
  The constant factor $n$ in the objective function is not relevant here, but it does make the dual linear program \eqref{lp:dual} studied below look more natural.
  If $(t^*,x^*)$ is an optimal solution of \eqref{lp:prim}, then $x^*\in\FW(V)$, and $t_i^*=\tfrac{1}{n}d_\triangle(v_i,x^*)$.
  Conversely, each Fermat--Weber point arises in this way.
  The constraints $v_{ij} \leq t_i + x_j$ describe the $\max$-tropical version of the envelope $\cE(V)$; see \cite[\textsection 6.1]{ETC}.
  In that reference that inequality would look like \enquote{$v_{ij} \leq -t_i + x_j$}.
  Yet the linear substitution $t_i\mapsto -t_i$ is natural here, because \eqref{lp:prim} is a minimization problem; see also \cite[Remark 6.28]{ETC}.
  The cells of $\CovDec(V)$ are precisely the projections of the faces of the unbounded ordinary polyhedron $\cE(V)\subset\RR^m\times\RR^n$ onto the $x$-coordinates; see \cite[Proposition 6.11]{ETC}.
  Let $F$ be the optimal face of the linear program \eqref{lp:prim}.
  The set $\FW(V)$ is the projection of $F$.

  The face $F$ is gotten by minimizing $t_{1} + \dots +  t_{m}$, or equivalently $t_{1} + \dots +  t_{m} + x_1 + \dots + x_n$, as the points $x$ are restricted to the hyperplane $\cH$.
  Let $D$ be the cell in the triangulation of $\Delta_{m-1}\times\Delta_{n-1}$ which is dual to $F$.
  Then $D$ contains the vector $(\tfrac{1}{m}\1,\tfrac{1}{n}\1)\in\RR^m\times\RR^n$ in its relative interior.
  In other words, $\FW(V)$ is dual to the central cell $C_V$.

  A sublevel set of a function $f:\torus{n}\to\RR$ is a set of the form $\SetOf{x\in\torus{n}}{f(x)\leq\alpha}$ for some $\alpha\in\RR$.
  If all of its sublevel sets are bounded, then the set of minima is bounded.
  We use this property to show that $\FW(V)$ is bounded.
  
  The sublevel sets of $d_{\triangle}(v_i,\cdot)$ are simplices if non-empty; in particular, they are bounded.
  Consequently, the function $x\mapsto\sum_{i\in[m]}d_{\triangle}(v_i,x)$ has bounded sublevel sets.
  The latter implies that $\FW(V)$ is bounded, as it is the set of minimizers of the aforementioned function.
  The $\max$-tropical convex hull of the rows of $V$ equals the union of the bounded covector cells in $\CovDec(V)$; see \cite[Corollary 6.17]{ETC}.
\end{proof}

\begin{example}\label{exmp:FW-prim}
  For the matrix $V$ from Example~\ref{exmp:FWpoint} the unique optimal solution of the primal linear program \eqref{lp:prim} reads
  \[
    t^*=(5,4,5,7,3) \quad \text{and} \quad x^*=(9,-6,-3) \enspace ,
  \]
  with optimal value $3 \cdot 24=72$.
  We have $\FW(V)=\{x^*\}$; see Fig.~\ref{fig:FWpoint}.
  As $\gcd(5,3)=1$, the uniqueness is implied by Theorem~\ref{thm:dimFW} below.
  The point $x^*$, which is a pseudovertex of $\CovDec(V)$, is dual to the central cell $C_V=\conv\{113,122,212,221\}$; see Fig.~\ref{fig:dualSubd}.
\end{example}

\begin{remark}
  Theorem~\ref{thm:convFW} reveals a similarity to the Euclidean case:
  by \cite[Proposition~19.1]{BMS:1999} any Fermat--Weber point with respect to the Euclidean distance is contained in the convex hull of the sites.
  The analogous result to Theorem~\ref{thm:convFW} for the symmetric tropical distance is \cite[Proposition 26]{ConvTreeSp}.
  As shown in \cite[Example 27]{ConvTreeSp}, in general, the symmetric tropical Fermat--Weber set is not contained in the tropical convex hull.
\end{remark}

Via the Cayley trick, the central cell $C_V$ in $\Sigma(V)$ corresponds to the central covector cell in $\CovDec(V)$, which is $\FW(V)$.
The dimension of the latter is the codimension of the former.

\begin{theorem} \label{thm:dimFW}
  The dimension of $\FW(V)$ is at most $\gcd(m,n)-1$.
  In particular, if $m$ and $n$ are relatively prime, the Fermat--Weber point is unique.
\end{theorem}

\begin{proof}
  We consider the regular subdivision, $\Sigma(V)$, of $\Delta_{m-1}\times\Delta_{n-1}$ induced by the matrix $V$.
  Let $C_V$ be the central cell.
  By \cite[\textsection 6.2]{DLRS} the vertices of $C_V$ are in bijection with the edges in a subgraph $G(C_V)$ of the complete bipartite graph $K_{m,n}$ with $c:=\codim(C_V)+1$ connected components; see also \cite[\textsection 4.7]{ETC}.
  We will show that $c\leq\gcd(m,n)$.
  Here we may assume that $V$ is generic, whence $C_V$ is a simplex; note that any refinement of $\Sigma(V)$ can only increase the codimension of the cell containing a specific point in its relative interior.
  If $c=1$, there is nothing to prove.
  Hence, we assume that $c\geq 2$.
	
  We consider a maximal simplex, $U$, of $\Sigma(V)$ which contains $C_V$, and let $T$ be the subtree of $K_{m,n}$ corresponding to $U$.
  A facet, $F$, of $U$ corresponds to the removal of some edge $e$ of $T$, yielding disjoint unions $[m]=I'\sqcup I''$ and $[n]=J'\sqcup J''$ such that the connected components of $T\setminus e$ are the restrictions on $I'\sqcup J'$ and $I''\sqcup J''$ and $e$ is incident to a vertex in $I''$ and one in $J'$.
  The hyperplane $\sum_{i'\in I'}x_{i'}+\sum_{j''\in J''}y_{j''}=1$ contains the facet $F$.
  If $I'$ were empty, then $\sum_{j'\in J'}y_{j'}=0$ for every $(x,y)\in C_V$.
  But $C_V$ contains the point $(\tfrac{1}{m}\1,\tfrac{1}{n}\1)$, so $J'$ must be empty as well.
  This contradicts the fact that $J'$ contains a vertex incident to $e$.
  Hence, $I'\neq\emptyset$ and, similarly, $J''\neq\emptyset$.

  Consequently, there exist two proper partitions $[m]=I_1\sqcup I_2\sqcup\dots\sqcup I_c$ and $[n]=J_1\sqcup J_2\sqcup\dots\sqcup J_c$ such that the restriction of $G(C_V)$ on $I_i\sqcup J_i$ is a tree, and there exist edges $e_2,\dots,e_c$ in $K_{m,n}$ between a vertex of $I_{i}$ and $J_{i-1}$ such that adding these edges we obtain the tree $T$.

  A supporting hyperplane for the facet corresponding to $T\setminus\{e_i\}$ is given by the equation
  \[\sum_{i\in I_1\sqcup\dots\sqcup I_{i-1}}x_i+\sum_{j\in J_i\sqcup\dots\sqcup J_c}y_j \ = \ 1 \enspace ,\]
  where $2\leq i\leq c$.
  Now the point $(\tfrac{1}{m}\1,\tfrac{1}{n}\1)$ is contained in these hyperplanes, yielding the $c-1$ equalities
  \[\tfrac{1}{m}\sum_{k<i}|I_k|+\tfrac{1}{n}\sum_{\ell\geq i}|J_\ell| \ = \ 1 \,, \qquad \text{for } 2\leq i\leq c \enspace .\]
  Multiplying with the least common multiple of $m$ and $n$ we obtain
  \[\tfrac{n}{\gcd(m,n)}\sum_{k<i}|I_k|+\tfrac{m}{\gcd(m,n)}\sum_{\ell\geq i}|J_\ell| \ = \ \lcm(m,n) \,, \qquad \text{for } 2\leq i\leq c \enspace .\]
  Putting those relations together with $\sum_{k\in [c]}|I_k|=m$ and $\sum_{\ell\in[c]}|J_\ell|=n$, we yield that $|I_k|$ is a multiple of ${m}/{\gcd(m,n)}$ for every $k\in[c]$ and $|J_\ell|$ is a multiple of ${n}/{\gcd(m,n)}$ for all $\ell\in[c]$.
  Further $I_k\neq\emptyset$ for every $k\in[c]$, and so $|I_k|\geq{m}/{\gcd(m,n)}$ for all $k\in[c]$.
  Hence
  \[ m \ = \ \sum_{k\in[c]}|I_k| \ \geq \ c\cdot\tfrac{m}{\gcd(m,n)} \enspace ,\]
  which implies $c\leq\gcd(m,n)$ as claimed.
\end{proof}

Restricting to the one-dimensional case (i.e., $n=2$), we recover the known fact that an odd number of points have a unique median, while the median can be selected from an interval for an even number of points.
The following example shows that the inequality in Theorem \ref{thm:dimFW} is tight for all $m$ and $n$; see also Example~\ref{exmp:FWpoint}.

\begin{example}\label{exmp:staircase}
  Our example employs the matrix $V=(v_{ij})\in\RR^{m\times n}$ with $v_{ij}=(i-1)(j-1)$.
  The rows are points on the tropical moment curve, and their ($\max$-)tropical convex hull is a tropical cyclic polytope; see \cite[\textsection 4]{BlockYu:2006} and \cite[Example~5.18]{ETC}.
  The dual triangulation $\Sigma(V)$ of $\Delta_{m-1}\times\Delta_{n-1}$ is the \emph{staircase triangulation}; see \cite[\textsection 6.2.3]{DLRS}.
  We explain the construction.
  
  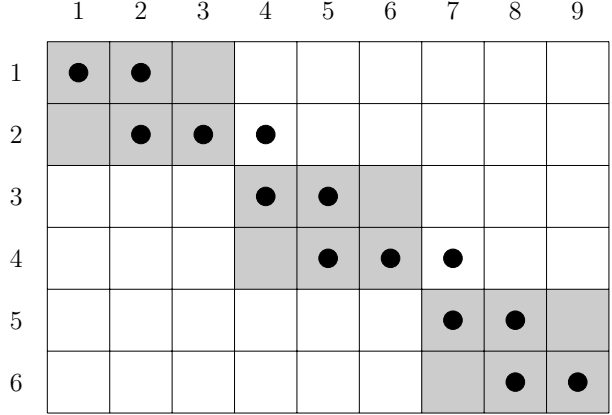
\begin{figure}[th]  \centering
		\resizebox{8cm}{!}{\begin{tikzpicture}

\fill[gray!40] (0,6) rectangle (3,4);
\fill[gray!40] (3,4) rectangle (6,2);
\fill[gray!40] (6,2) rectangle (9,0);

\draw (0,0) grid (9,6);

\node at (0.5,5.5) [circle, fill = black, inner sep = 3pt, draw = black] {};
\node at (1.5,5.5) [circle, fill = black, inner sep = 3pt, draw = black] {};
\node at (1.5,4.5) [circle, fill = black, inner sep = 3pt, draw = black] {};
\node at (2.5,4.5) [circle, fill = black, inner sep = 3pt, draw = black] {};

\node at (3.5,4.5) [circle, fill = black, inner sep = 3pt, draw = black] {};

\node at (3.5,3.5) [circle, fill = black, inner sep = 3pt, draw = black] {};
\node at (4.5,3.5) [circle, fill = black, inner sep = 3pt, draw = black] {};
\node at (4.5,2.5) [circle, fill = black, inner sep = 3pt, draw = black] {};
\node at (5.5,2.5) [circle, fill = black, inner sep = 3pt, draw = black] {};

\node at (6.5,2.5) [circle, fill = black, inner sep = 3pt, draw = black] {};

\node at (6.5,1.5) [circle, fill = black, inner sep = 3pt, draw = black] {};
\node at (7.5,1.5) [circle, fill = black, inner sep = 3pt, draw = black] {};
\node at (7.5,0.5) [circle, fill = black, inner sep = 3pt, draw = black] {};
\node at (8.5,0.5) [circle, fill = black, inner sep = 3pt, draw = black] {};

\node at (0.5,6.5) [] {$1$};
\node at (1.5,6.5) [] {$2$};
\node at (2.5,6.5) [] {$3$};
\node at (3.5,6.5) [] {$4$};
\node at (4.5,6.5) [] {$5$};
\node at (5.5,6.5) [] {$6$};
\node at (6.5,6.5) [] {$7$};
\node at (7.5,6.5) [] {$8$};
\node at (8.5,6.5) [] {$9$};

\node at (-0.5,5.5) [] {$1$};
\node at (-0.5,4.5) [] {$2$};
\node at (-0.5,3.5) [] {$3$};
\node at (-0.5,2.5) [] {$4$};
\node at (-0.5,1.5) [] {$5$};
\node at (-0.5,0.5) [] {$6$};

\end{tikzpicture}}
    \caption{%
      Staircase in a $6\times 9$ grid}
    \label{fig:staircase}
  \end{figure}

  In this case, we represent the simplices as points in an $m{\times}n$ grid instead of subgraphs of $K_{m,n}$.
  The staircase triangulation consists of all the paths in an $m{\times}n$ grid starting at $(1,1)$ and ending at $(m,n)$ obtained by going right or down.
  Figure \ref{fig:staircase} portrays such a path.
  Due to its shape, it is called a \emph{staircase}.
  We will show that in the staircase triangulation the simplex containing $(\tfrac{1}{m}\1,\tfrac{1}{n}\1)$ lies in a cell of codimension $\gcd(m,n)-1$.
  For improved readability, we abbreviate $d=\gcd(m,n)$, $a={m}/{d}$, and $b={n}/{d}$.
  
  If $d=1$, then there is a unique maximal simplex in the staircase triangulation containing $(\tfrac{1}{m}\1,\tfrac{1}{n}\1)$ in its interior.
  To find the precise staircase, we use the Northwest Corner Rule \cite[\textsection 8.3.1]{DantzigThapa:1997} in the standard transportation array with marginal column $\tfrac{1}{m}\1$ and marginal row $\tfrac{1}{n}\1$.
  The visited cells form the staircase, which we call the \emph{central staircase}.
  Note that this method provides also barycentric coordinates for $(\tfrac{1}{m}\1,\tfrac{1}{n}\1)$ in this simplex.
  
  When $d\geq 2$, consider the partitions in $d$ subsets $[m]=I_1\sqcup\dots\sqcup I_d$ and $[n]=J_1\sqcup\dots\sqcup J_d$ where $I_i=\{(i-1)a+1,(i-1)a+2,\dots,ia\}$ and $J_i=\{(i-1)b+1,(i-1)b+2,\dots,ib\}$.
  Consider on $I_i\times J_i$ the central staircase on an $a\times b$ grid and add the points $(ia,ib+1)$ for $i=1,\dots,d-1$: in Figure \ref{fig:staircase} the blocks on $I_i\times J_i$ correspond to the gray areas whereas the added points are those in the white area.
  This staircase corresponds to a maximal simplex $U$ in the staircase triangulation, which contains $(\tfrac{1}{m}\1,\tfrac{1}{n}\1)$.
  
  The removal of the grid point $(ia,ib+1)$ yields the facet defined by
  \[ \sum_{k\leq ia}x_k+\sum_{\ell>ib}y_\ell \ = \ 1 \enspace. \]
  Using $d={m}/{a}={n}/{b}$, we see that $(\tfrac{1}{m}\1,\tfrac{1}{n}\1)$ is contained in this facet.
  In total, there are $d-1$ facets of this form.
  
  Each remaining facet of $U$ corresponds to the removal of a grid point from some block $I_k\times J_k$.
  In particular, each facet induces a partition $[m]=I'\sqcup I''$ such that $I_1\sqcup\dots\sqcup I_{k-1}\subset I'$, $I_{k+1}\sqcup\dots\sqcup I_{d}\subset I''$ for some $k\in[d]$; similarly, there is a partition $J'\sqcup J''$ on $[n]$.
  Moreover, at least one of the intersections $I_k\cap I'\cap I''$ and $J_k\cap J'\cap J''$ is nonempty, which implies that at least one of $|I'|$ and $|J''|$ is not a multiple of $d$.
  As in the proof of Theorem~\ref{thm:dimFW}, the facet defining equation
  \[ \sum_{i'\in I'}x_{i'}+\sum_{j''\in J''}y_{j''} \ = \ 1 \]
  is not satisfied by $(\tfrac{1}{m}\1,\tfrac{1}{n}\1)$.
  Since $(\tfrac{1}{m}\1,\tfrac{1}{n}\1)$ is contained in precisely $d-1$ facets of $U$, the dimension of $\FW(V)$ equals $d-1$.
  
  The staircase for $U$ is obtained by using the Northwest Corner Rule with breaking ties by going East.
  If we break the ties randomly, then $2^{d-1}$ staircases appear with nonzero probability.
  These staircases are in bijection with the ordinary vertices of $\FW(V)$, which is a $(d{-}1)$-dimensional cube, seen as an ordinary polytope.
\end{example}


\begin{remark}
  In the special case $m=n-1$ computing tropical Fermat--Weber points reduces to the tropical Cramer rule \cite[\textsection 4.9]{ETC}. 
  Adding one more site, $p$, the new Fermat--Weber set, $F$, can have  
  higher dimension, but contains the tropical Cramer point, $c$.  
  Perturbing $c$ in the direction of $p$, we obtain a point in the relative interior of $F$. This also agrees with results of G\"artner and Jaggi \cite[\textsection 4.1]{GaertnerJaggi:2006} in the context of \enquote{tropical support vector machines} on computing a \enquote{separating hyperplane for $n$ points}.
\end{remark}

\begin{remark}
  As a consequence of \cite[Theorem~6.14]{ETC} the Fermat--Weber set $\FW(\transpose{V})$ in $\torus{m}$ of the $n$ columns of $V$ is affinely isomorphic to the Fermat--Weber set $\FW(V)$ in $\torus{n}$ of the $m$ rows.
\end{remark}

\begin{remark}
  Our analysis rests on the decision, in \eqref{eq:from-the-sites}, to look at the distances from the sites to the Fermat--Weber points.
  This leads to $\min$-tropical hyperplane arrangements and $\max$-tropical convex hulls.
  Reversing the direction, i.e., considering distances from the Fermat--Weber points to the sites, amounts to exchanging $\min$ and $\max$ throughout.
  Conceptually, the results remain the same; cf.\ \cite[\textsection 1.3]{ETC}.
\end{remark}

\section{Transportation} \label{sec:transportation}
\noindent
This section comprises the algorithmic core of this paper.
The basic ingredient is the transportation problem which already occurred in Example~\ref{exmp:staircase}.

Again let us fix a matrix $V=(v_{ij})\in\RR^{m\times n}$.
Whenever it will suit us, we may also view the rows of $V$ is a $m$ labeled points $v_1,\dots,v_m$ in $\cH$ or $\torus{n}$.
The following is the dual of the linear program \eqref{lp:prim} with variables $\lambda$ and $y_{ij}$ for $i\in[m]$ and $j\in[n]$:
\begin{equation} \label{lp:dual}
  \begin{array}{ll@{}ll}
    \text{maximize}  & \displaystyle \sum_{i\in[m]}\sum_{j\in[n]}v_{ij}\cdot y_{ij} & \\
    \text{subject to}& \displaystyle \sum_{j\in[n]} y_{ij} \ = \ n \enspace , & \quad \text{for } i\in[m] \\
                     & \displaystyle \lambda + \sum_{i\in[m]} y_{ij} \ = \ 0 \enspace , & \quad \text{for } j\in[n] \\
                     & \displaystyle y_{ij} \ \geq \ 0 \enspace , & \quad \text{for } i\in[m] \text{ and } j\in[n] \enspace .
  \end{array}
\end{equation}
From $0 = n\cdot \lambda + \sum_{i,j}y_{ij} = n\cdot (\lambda + m)$, we get $\lambda=-m$ for every feasible point.
By substituting that value in \eqref{lp:dual} we obtain the standard linear programming formulation of a transportation problem; see \cite[\textsection 21.6]{Schrijver03:CO_A}.
The primal linear program \eqref{lp:prim} and the envelope $\cE(V)$ are related to transportation via Hitchcock's theorem \cite[Theorem 21.13]{Schrijver03:CO_A}.

As the right hand sides are the integral constants $m$ and $n$, the feasible region of \eqref{lp:dual} is a \emph{central transportation polytope}, and we denote it $T(m,n)$.
The polytope $T(m,n)$ is known to be integral \cite[Lemma 2.13]{DeLoera+Kim:2014}.
The nonzero entries of any vertex, which is an $m{\times}n$-matrix, defines a subgraph of $K_{m,n}$ by picking edges \cite[Lemma 2.9]{DeLoera+Kim:2014}.
This is the \emph{support graph} of that vertex, and this is a forest; see \cite[Theorem 21.15]{Schrijver03:CO_A}.

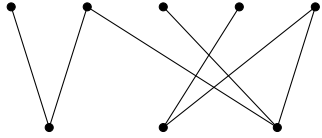
\begin{figure}[th] \centering
  \begin{tikzpicture}
    \newcommand\yhi{0.8}
    \newcommand\ylo{-0.8}
    \newcommand\sz{3pt}
    \coordinate (r1) at (-2,\yhi);
    \coordinate (r2) at (-1,\yhi);
    \coordinate (r3) at (0,\yhi);
    \coordinate (r4) at (1,\yhi);
    \coordinate (r5) at (2,\yhi);
    \coordinate (c1) at (-1.5,\ylo);
    \coordinate (c2) at (0,\ylo);
    \coordinate (c3) at (1.5,\ylo);
    \
    \foreach \i/\j in {1/1,2/1,2/3,3/3,4/2,5/2,5/3}{
      \draw (r\i) -- (c\j);
    }
    \foreach \i in {1,2,3,4,5}{
      \filldraw[blackdot] (r\i) circle (\sz);
    }
    \foreach \j in {1,2,3}{
      \filldraw[blackdot] (c\j) circle (\sz);
    }
  \end{tikzpicture}
  \caption{Spanning subtree of $K_{5,3}$ encoding the covector of the unique Fermat--Weber point from Example~\ref{exmp:FWpoint}.  Column nodes at the bottom}
  \label{fig:covector}
\end{figure}

\begin{example}\label{exmp:covector}
  For the $5{\times}3$-matrix $V$ from Example~\ref{exmp:FWpoint} the transpose of the unique optimal solution of \eqref{lp:dual} is
  \begin{equation}\label{eq:covector}
    \transpose{(y_{ij}^*)} \ = \
    \begin{pmatrix}
      3 & 2 & 0 & 0 & 0 \\
      0 & 0 & 0 & 3 & 2 \\
      0 & 1 & 3 & 0 & 1
    \end{pmatrix}
    \enspace .
  \end{equation}
  The support graph a spanning tree; see Fig.~\ref{fig:covector}.
  By Theorem~\ref{thm:FW-dual} below that tree encodes the covector of the unique Fermat--Weber point $x^*$.
  The degree sequence for the column nodes is $(223)$.
  This is the coarse type of $x^*$ and also the componentwise maximum of the four vertices of the central cell in the mixed subdivision; cf.\ Example~\ref{exmp:FW-prim} and \cite[\textsection 4.5]{ETC}.
\end{example}

\begin{figure}[th]
  \centering
	\resizebox{8cm}{!}{\begin{tikzpicture}[scale=1.5,ultra thick
,x={(0.866,-0.5)},y={(0,1)},z={(-0.866,-0.5)}]
		\definecolor{tolorange}{RGB}{238,119,51} 
		\definecolor{tolblue}{RGB}{0,119,187}    
		\definecolor{tolgreen}{RGB}{0,153,136}   
		\definecolor{tolpurple}{RGB}{136,34,85}  
		
		\draw[dotted] (-5,-9) -- (-5,-4) -- (-10,-4);
		\draw[dotted] (-5,-4) -- (0,1);
		
		\node at (-2,-1) [circle, fill = tolgreen, draw = white, semithick, inner sep = 2.5pt, label={[tolgreen]right:$(3,-3,0)$}] {};
		\node at (-9,-4) [circle, fill = tolorange, draw = white, semithick, inner sep = 2.5pt, label={[tolorange]right:$(13,-14,1)$}] {};
		\node at (-8,-3) [circle, fill = tolpurple, draw = white, semithick, inner sep = 2.5pt, label={[tolpurple]right:$(11,-13,2)$}] {};
		\node at (-3,-7) [circle, fill = black, draw = white, semithick, inner sep = 2.5pt, label={[black]right:$(10,1,-11)$}] {};
		\node at (-7,-7) [circle, fill = tolblue, draw = white, semithick, inner sep = 2.5pt, label={[tolblue]right:$(14,-7,-7)$}] {};
		
		\node (FW) at (-5,-4) [circle, draw= black, fill = white, semithick, inner sep = 2.5pt, label={[black]right:$(9,-6,-3)$}] {};
\end{tikzpicture}}
	
  \caption{Five sites in $\torus{3}$ and the unique point that evenly splits them.
    Each closed sector of the max-tropical halfspace with apex $(9,-6,-3)$ contains at least two sites}
  \label{fig:evenSplit}
\end{figure}
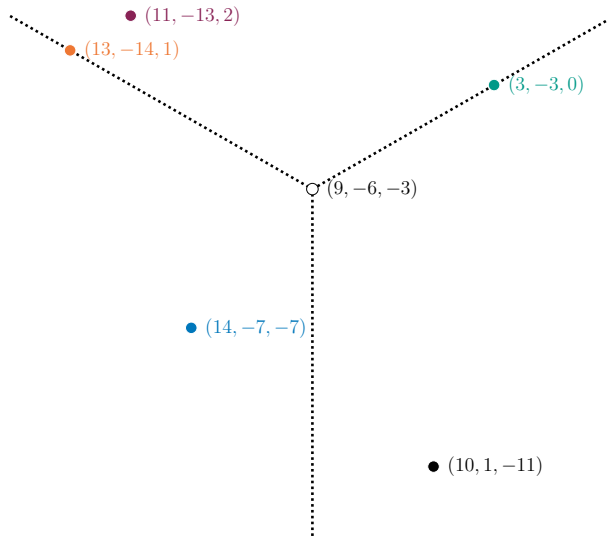

Assuming $m\geq n$, Tokuyama and Nakano gave an algorithm (which they called \enquote{splitter finding}) to solve a transportation problem like \eqref{lp:dual} in $O(n^2m\log^2 m)$ time \cite[Theorem~3.1]{Tokuyama+Nakano:1995}.
We borrow some of their ideas.
For $J\subseteq[n]$ and $u\in\cH$ consider the set
\[ S_J(u) \ := \ \bigSetOf{x\in\cH}{\max_{j\in J}(x_j-u_j)\geq\max_{i\notin J}(x_i-u_i)} \enspace . \]
We have $S_\emptyset(u)=\emptyset$, and $S_{[n]}(u)=\cH$, where we use the convention that the maximum of the empty set is $-\infty$, the neutral element of the tropical addition.
If $J$ is a nonempty proper subset of $[n]$, then $S_J(u)$ is a max-tropical halfspace with apex $u$; see \cite[\textsection 7.1]{ETC}.
In the special case where $J=\{j\}$ is a singleton that tropical halfspace is a closed sector; in general, $S_J(u)=\bigcup_{j\in J}S_{\{j\}}(u)$.
\begin{definition}
  We say that $u\in\cH$ \emph{evenly splits} $V$ if for every subset $J$ of $[n]$ we have $n\cdot|V\cap S_J(u)|\geq m\cdot |J|$.
\end{definition}
Tokuyama and Nakano \cite{Tokuyama+Nakano:1995} call the point $u$ a \enquote{$\1$-splitter}, and the sectors are the \enquote{regions split by~$u$}.
Theorem \ref{thm:FW-dual} below may be seen as our interpretation of their results in the setting of tropical convexity.

\begin{example} \label{exmp:evenSplit}
  Consider the points $V$ from Example \ref{exmp:FWpoint}.
  The Fermat--Weber point $u=(9,-6,-3)$ evenly splits them.
  This is illustrated in Figure \ref{fig:evenSplit}, where we have drawn the $\max$-tropical hyperplane based at $u$ with dotted lines.
  In particular, we see the subdivision of $\torus{3}$ in three convex regions.
\end{example}

\begin{theorem}\label{thm:FW-dual}
  A point $u\in\cH$ evenly splits $V$ if and only if $u\in\FW(V)$.
  The support graph of an optimal dual solution $y^*$ is the covector of the Fermat--Weber point $u(y^*)$.
\end{theorem}

\begin{proof}
  Assume that $u$ evenly splits $V$ and consider $t_i=\max_{j\in[n]}(v_{ij}-u_j)$ for $i\in[m]$.
  Thus $(u,t)$ is a feasible solution of the primal linear program (\ref{lp:prim}).
  Moreover, $t_i+u_j=v_{ij}$ if and only if $v_i\in S_{\{j\}}(u)$.
  
  Now \cite[Theorem 2.2]{Tokuyama+Nakano:1995} says that there exists a solution $(y_{ij})$ of the transportation problem (\ref{lp:dual}), which is dual to (\ref{lp:prim}), such that $(u,t)$ and $(y_{ij})$ satisfy the complementary slackness condition.
  Indeed, if $t_i+u_j\neq v_{ij}$, then $v_i\notin S_{\{j\}}(u)$, so the aforementioned result gives $y_{ij}=0$.
  So it follows from linear programming duality that $(u,t)$ is an optimal solution of (\ref{lp:prim}).
  In particular, $u$ is a Fermat--Weber point for $V$.
  
  For the converse, we denote by $\phi$ the convex function $\tfrac{1}{n}\sum_{s\in V}d_\triangle(s,\cdot)$.
  Also, for every subset $J$ of $[n]$, denote by $\sigma_J$ the cardinality of $V\cap S_J(u)$.
  Abbreviating $f_J:=(n-|J|)\sum_{j\in J}e_j-|J|\sum_{i\in[n]\setminus J}e_i$, which is a point in $\cH$, we obtain:
  \begin{itemize}
  \item if $s\in S_J(u)$, then $s\in S_J(u-\delta f_J)$ for any $\delta\geq 0$;
  \item if $s\notin S_J(u)$, then $s\notin S_J(u-\delta f_J)$ for any $\delta\geq 0$ sufficiently small.
  \end{itemize}
  The condition $s\in S_J(u)$ implies $d_\triangle(s,u)=n(s_j-u_j)$ for some $j\in J$.
  Therefore, for $\delta>0$ sufficiently small, we obtain
  \begin{itemize}
  \item $d_\triangle(s,u-\delta f_J)=d_\triangle(s,u)+n\cdot\delta(n-|J|)$, when $s\in S_J(u)$;
  \item $d_\triangle(s,u-\delta f_J)=d_\triangle(s,u)-n\cdot\delta|J|$, when $s\notin S_J(u)$.
  \end{itemize}
  Summing up and dividing by $n$ yields
  \begin{equation} \label{eq:phi_fJ}
    \phi(u-\delta f_J) \ = \ \phi(u)+\delta(n-|J|)\sigma_J-\delta|J|(m-\sigma_J) \ = \ \phi(u)+\delta\left(n\sigma_J-m|J|\right)
  \end{equation}
  for any $\delta> 0$ sufficiently small.
  If $u\in\FW(V)$, then $u$ is a minimizer, so $\phi(u-\delta f_J)\geq \phi(u)$ for any $\delta>0$ and $J\subseteq [n]$.
  Equation (\ref{eq:phi_fJ}) implies $n\sigma_J\geq m|J|$ for every $J\subseteq[n]$, under this assumption.
  Consequently, $u$ evenly splits $V$.
\end{proof}

The argument in the proof of Theorem \ref{thm:FW-dual} leads to the following algorithm for obtaining a tropical Fermat--Weber point from an optimal solution $y^*$ of  (\ref{lp:dual}).
By complementary slackness each edge $(i,j)$ of the support graph, $B(y^*)$, imposes the equality
\begin{equation} \label{eq:slackness}
t_i+x_j \ = \ v_{ij}
\end{equation}
for any optimal dual solution.
If we assume $x_j=0$ for some column node $j$, then we can perform a depth-first search from $j$ and recover all the the other values of $(t,x)$ using the equalities in (\ref{eq:slackness}).
There may be more connected components, in which case we start a depth-first search from every unvisited column node.
In this way all the values are recovered, as $B(y^*)$ is spanning---no row or column of $y^*$ can be zero as its elements must sum up to a positive number.
Note that the point $x$ obtained this way may not lie in $\cH$.
Yet, by adding $(x_1+\dots+x_n)/n$ to every entry of $x$ and subtracting the same value from every entry of $t$, we obtain a feasible solution $(t^*,x^*)$ of (\ref{lp:prim}) that satisfies the equations (\ref{eq:slackness}).
In particular, that solution is optimal, whence $x^*\in\FW(V)$.
So the method of \cite{Tokuyama+Nakano:1995} gives the following complexity result.
\begin{corollary} \label{cor:compFW}
  Assuming $m\geq n$, one tropical Fermat--Weber point can be computed in $O(n^2m\log^2 m)$ time.
\end{corollary}
The complexity bounds of the algorithms by Kleinschmidt and Schannath \cite{Kleinschmidt+Schannath:1995} and Brenner \cite{Brenner:2008} are similar, but slightly different.
Naturally, they carry over as well.

It is also natural to ask for explicit representations of the entire Fermat--Weber set.
A first idea could be to list the vertices of $\FW(V)$, seen as an ordinary polytope.
Example~\ref{exmp:staircase} shows that cubes occur as Fermat--Weber sets, and their number of vertices is exponential in the dimension.
Yet, the polytropal structure allows for more efficient choices.
Namely, a polytope in $\torus{n}$ has at most $n$ tropical vertices and at most $n^2-n$ ordinary facets.

Let us start with the latter.
Since we know the possible directions of the (outer) facet normal vectors, $e_k-e_\ell$ for $k\neq\ell$, we can find the ordinary facets by solving $n^2-n$ linear programs like:
\begin{equation} \label{eq:facetFW}
  \begin{array}{ll@{}ll}
    \text{maximize}  & \displaystyle  x_k - x_\ell& \\
    \text{subject to}& \displaystyle v_{ij} - x_j \leq t_i \,, & \quad \text{for } i\in[m] \text{ and } j\in[n] \\
                     & \displaystyle x_{1} + \dots + x_{n} = 0  & \\
                     & \displaystyle t_{1} + \dots + t_{m} = p^*/n &
    \end{array}
\end{equation}
where $p^*$ is the optimal value of (\ref{lp:prim}).
The constraint matrix has only $0$ and $\pm 1$ entries, and so a linear program of the form (\ref{eq:facetFW}) can be solved in strongly polynomial time \cite[Corollary~15.3a]{Schrijver:1986}.

Each such linear program yields one tight inequality $x_k-x_\ell\leq a_{k,\ell}$ of $\FW(V)$, where $a_{k,\ell}$ is the optimal value.
Then, the tropical vertices can be found as $A_k=(-a_{k,1},\dots,-a_{k,n})\in\torus{n}$, where $a_{k,k}=0$.
From Corollary \ref{cor:compFW} we thus infer the following result.
	
\begin{corollary} \label{cor:strongPoly}
  The ordinary facet description and the tropical vertices of $\FW(V)$ can be found in strongly polynomial time.
\end{corollary}
If the tropical vertices are known, then we can check if a given point lies in $\FW(V)$ in $O(n^2)$ time by checking the criterion \cite[Proposition 5.37]{ETC}.

\section{Tropical median consensus trees}
\label{sec:trees}
\noindent
Phylogenetics is a branch of (computational) biology that seeks to associate trees to mark ancestral relations among given \emph{taxa}, e.g., species \cite{SempleSteel:2003}.
The taxa correspond to the leaves.
Here we show how asymmetric tropical Fermat--Weber sets give rise to a new algorithm for finding consensus trees.
Since there is no particular shortage on consensus tree methods, we compare its features to other approaches, and we discuss the practical applicability.

\subsection{Equidistant trees}
We recall known facts about ultrametrics and equidistant trees to fix our notation; see \cite[\S7.2]{SempleSteel:2003} and \cite[\S10.9]{ETC}.
A \emph{dissimilarity map} is a symmetric $n{\times}n$-matrix $D=(d_{ij})$ with zero diagonal.
It is called an \emph{ultrametric} if $D$ is nonnegative, and the \emph{ultrametric inequality}
\begin{equation}\label{eq:ultrametric}
  d_{ik} \ \leq \ \max (d_{ij} , d_{jk})
\end{equation}
holds for all $i,j,k\in[n]$.
Since the (zero) diagonal is implicit, and the matrix is required to be symmetric, we may view a dissimilarity map as an element of $\RR^{\tbinom{n}{2}}$.

A rooted metric tree with $n$ labeled leaves is \emph{equidistant} if the distance from any leaf to the root is the same.
It is known that the ultrametrics are precisely the distance functions among the leaves in an equidistant tree; see \cite[Theorem 7.2.5]{SempleSteel:2003}.
Note that the all-ones vector $\1$ of length $\tbinom{n}{2}$ is an ultrametric.
The corresponding equidistant tree has interior edges of length zero, while each leaf edge has length $\tfrac{1}{2}$.

Ardila and Klivans \cite[Theorem 3]{ArdilaKlivans:2006} showed that a dissimilarity map is an ultrametric if and only if it corresponds to a point on the Bergman fan $\BergmanFan(K_n)$ of the complete graph $K_n$.
The Bergman fan of a matroid is a special case of a tropical linear space, i.e., with constant coefficients.
The ultrametric inequality \eqref{eq:ultrametric} is a tropical convexity condition with respect to $\max$.
From a dissimilarity map $D$ and any real constant $c$ we can construct the dissimilarity map $D+c\1$.
Moreover, if $D$ is an ultrametric, and $D+c\1$ is nonnegative, then it is an ultrametric, too.
In this way, we may view $\BergmanFan(K_n)$ as a $\max$-tropical linear space in the tropical projective torus $\torus{\tbinom{n}{2}}$; see \cite[Proposition 16]{ConvTreeSp} and \cite[Theorem~3]{YoshidaZhangZhang:2019}.
We abbreviate $\cT_n:=\BergmanFan(K_n)/\RR\1$ and call it the \emph{space of equidistant trees} on $n$ labeled leaves.
Now we can apply our results from Sections \ref{sec:fermat-weber} and \ref{sec:transportation} to points in $\cT_n$.
Our first observation says that Fermat--Weber points of equidistant trees are again equidistant trees.

\begin{theorem} \label{thm:trop_consensus}
  Let $V\subset\cT_n$ be a finite set of equidistant trees on $n$ leaves.
  Then the tropical polytope $\FW(V)$ is contained in $\cT_n$.
  Moreover, any two trees in $\FW(V)$ share the same tree topology.
\end{theorem}

\begin{proof}
  According to Theorem \ref{thm:convFW} the set $\FW(V)$ is contained in the $\max$-tropical convex hull of $V$.
  The space of equidistant trees $\cT_n$ is a $\max$-tropical linear space and thus $\max$-tropically convex; see \cite[Proposition 10.33]{ETC}.
  This is the first claim.
  Page, Yoshida and Zhang showed \cite[Theorem~3.2]{Page+Yoshida+Zhang:2020} that the trees in any relatively open cell of the covector decomposition of $V$ in $\cT_n$ share the same tree topology.
  With this observation the second claim follows also from Theorem~\ref{thm:convFW}.
\end{proof}

As we know, Fermat--Weber points are not unique, in general.
Here is a more precise statement.

\begin{corollary} \label{cor:trop_consensus}
  Let $V\subset\cT_n$ be a set of $m$ equidistant trees on $n$ leaves.
  Then
  \[
    \dim\FW(V) \ \leq \ \min\bigl(\, n-1 \,,\; \gcd(m,\tbinom{n}{2}) \,\bigr)-1 \enspace.
  \]
\end{corollary}

\begin{proof}
  This follows from Theorem \ref{thm:dimFW} and the fact that the graphic matroid of $K_n$ has rank $n-1$, so $\dim\cT_n= n-2$.
\end{proof}

\begin{figure}[th]\centering
    
    \begin{tabular}{cc}
      \resizebox{7.5cm}{!}{\begin{tikzpicture}[sloped]
    \node (A) at (0,0) {A};
    \node (B) at (1,0) {B};
    \node (C) at (2,0) {C};
    \node (D) at (3,0) {D};
    \node (E) at (4,0) {E};
    \node (F) at (5,0) {F};
    \node (G) at (6,0) {G};
    \node (H) at (7,0) {H};
    \node (I) at (8,0) {I};

    \node (CD) at (2.5,1) {};
    \node (BCD) at (1.75,2) {};
    \node (EF) at (4.5,1) {};
    \node (HI) at (7.5,1) {};
    \node (GHI) at (6.75,2) {};
    \node (EFGHI) at (5.625,4) {};
    \node (BCDEFGHI) at (3.6875,7) {};
    \node (all) at (1.84375,8) {};
    \node (root) at (1.84375,8.25) {};
    
    \draw  (A) |- (all.center);
    \draw  (B) |- (BCD.center);
    \draw  (C) |- (CD.center);
    \draw  (D) |- (CD.center);
    \draw  (E) |- (EF.center);
    \draw  (F) |- (EF.center);
    \draw  (G) |- (GHI.center);
    \draw  (H) |- (HI.center);
    \draw  (I) |- (HI.center);
    \draw  (CD.center) |- (BCD.center);
    \draw  (HI.center) |- (GHI.center);
    \draw  (EF.center) |- (EFGHI.center);
    \draw  (GHI.center) |- (EFGHI.center);
    \draw  (BCD.center) |- (BCDEFGHI.center);
    \draw  (EFGHI.center) |- (BCDEFGHI.center);
    \draw  (BCDEFGHI.center) |- (all.center);
    \draw  (all.center) |- (root);
    
    \draw (-1,0) -- (-1,8);
    
    \foreach \y in {0,1,2,3,4,5,6,7,8}                     
    \draw[shift={(0,\y)},color=black] (-1,0) -- (-1.1,0);
    
    \foreach \y in {0,1,2,3,4,5,6,7,8}   
    \node[left] at (-1.1,\y) {$\y$} ;
\end{tikzpicture}}
      & \resizebox{7.5cm}{!}{\begin{tikzpicture}[sloped]
    \node (A) at (0,0) {A};
    \node (B) at (2,0) {B};
    \node (C) at (1,0) {C};
    \node (D) at (3,0) {D};
    \node (E) at (4,0) {E};
    \node (F) at (5,0) {F};
    \node (G) at (6,0) {G};
    \node (H) at (7,0) {H};
    \node (I) at (8,0) {I};

    \node (BD) at (2.5,1) {};
    \node (CBD) at (1.75,2) {};
    \node (ACBD) at (0.875,3) {};
    \node (EF) at (4.5,1) {};
    \node (HI) at (7.5,1) {};
    \node (GHI) at (6.75,2) {};
    \node (EFGHI) at (5.625,4) {};
    \node (all) at (3.25,8) {};
    \node (root) at (3.25,8.25) {};
    
    \draw  (A) |- (ACBD.center);
    \draw  (B) |- (BD.center);
    \draw  (C) |- (CBD.center);
    \draw  (D) |- (BD.center);
    \draw  (E) |- (EF.center);
    \draw  (F) |- (EF.center);
    \draw  (G) |- (GHI.center);
    \draw  (H) |- (HI.center);
    \draw  (I) |- (HI.center);
    \draw  (BD.center) |- (CBD.center);
    \draw  (HI.center) |- (GHI.center);
    \draw  (EF.center) |- (EFGHI.center);
    \draw  (GHI.center) |- (EFGHI.center);
    \draw  (CBD.center) |- (ACBD.center);
    \draw  (ACBD.center) |- (all.center);
    \draw  (EFGHI.center) |- (all.center);
    \draw  (all.center) |- (root);
    
    \draw (-1,0) -- (-1,8);
    
    \foreach \y in {0,1,2,3,4,5,6,7,8}                     
    \draw[shift={(0,\y)},color=black] (-1,0) -- (-1.1,0);
    
    \foreach \y in {0,1,2,3,4,5,6,7,8}   
    \node[left] at (-1.1,\y) {$\y$} ;
\end{tikzpicture}} \\
      \small (a)  & \small (b) \\
      \resizebox{7.5cm}{!}{\begin{tikzpicture}[sloped]
    \node (A) at (0,0) {A};
    \node (B) at (5,0) {B};
    \node (C) at (1,0) {C};
    \node (D) at (2,0) {D};
    \node (E) at (3,0) {E};
    \node (F) at (4,0) {F};
    \node (G) at (6,0) {G};
    \node (H) at (7,0) {H};
    \node (I) at (8,0) {I};

    \node (CD) at (1.5,1) {};
    \node (ACD) at (0.75,2) {};
    \node (FB) at (4.5,1) {};
    \node (EFB) at (3.75,2) {};
    \node (GHI) at (7,2) {};
    \node (all) at (3.75,8) {};
    \node (root) at (3.75,8.25) {};
    
    \draw  (A) |- (ACD.center);
    \draw  (B) |- (FB.center);
    \draw  (C) |- (CD.center);
    \draw  (D) |- (CD.center);
    \draw  (E) |- (EFB.center);
    \draw  (F) |- (FB.center);
    \draw  (G) |- (GHI.center);
    \draw  (H) |- (GHI.center);
    \draw  (I) |- (GHI.center);
    \draw  (CD.center) |- (ACD.center);
    \draw  (FB.center) |- (EFB.center);
    \draw  (GHI.center) |- (all.center);
    \draw  (ACD.center) |- (all.center);
    \draw  (EFB.center) |- (all.center);
    \draw  (all.center) |- (root);
    
    \draw (-1,0) -- (-1,8);
    
    \foreach \y in {0,1,2,3,4,5,6,7,8}                     
    \draw[shift={(0,\y)},color=black] (-1,0) -- (-1.1,0);
    
    \foreach \y in {0,1,2,3,4,5,6,7,8}   
    \node[left] at (-1.1,\y) {$\y$} ;
	\end{tikzpicture}}
      & \resizebox{7.5cm}{!}{\begin{tikzpicture}[sloped]
    \node (A) at (0,0) {A};
    \node (B) at (1,0) {B};
    \node (C) at (2,0) {C};
    \node (D) at (3,0) {D};
    \node (E) at (4,0) {E};
    \node (F) at (5,0) {F};
    \node (G) at (6,0) {G};
    \node (H) at (7,0) {H};
    \node (I) at (8,0) {I};

    \node (CD) at (2.5,1) {};
    \node (EF) at (4.5,1) {};
    \node (HI) at (7.5,1) {};
    \node (GHI) at (6.75,2) {};
    \node (BCDEFGHI) at (3.6875,7) {};
    \node (all) at (1.84375,8) {};
    \node (root) at (1.84375,8.25) {};
    
    \draw  (A) |- (all.center);
    \draw  (B) |- (BCDEFGHI.center);
    \draw  (C) |- (CD.center);
    \draw  (D) |- (CD.center);
    \draw  (E) |- (EF.center);
    \draw  (F) |- (EF.center);
    \draw  (G) |- (GHI.center);
    \draw  (H) |- (HI.center);
    \draw  (I) |- (HI.center);
    \draw  (CD.center) |- (BCDEFGHI.center);
    \draw  (HI.center) |- (GHI.center);
    \draw  (EF.center) |- (BCDEFGHI.center);
    \draw  (GHI.center) |- (BCDEFGHI.center);
    \draw  (BCDEFGHI.center) |- (all.center);
    \draw  (all.center) |- (root);
    
    \draw (-1,0) -- (-1,8);
    
    \foreach \y in {0,1,2,3,4,5,6,7,8}                     
    \draw[shift={(0,\y)},color=black] (-1,0) -- (-1.1,0);
    
    \foreach \y in {0,1,2,3,4,5,6,7,8}   
    \node[left] at (-1.1,\y) {$\y$} ;
	\end{tikzpicture}}
 \\
      \small (c)  & \small (d) Consensus tree 
    \end{tabular}

    \caption{Three trees (a),(b),(c) and their tropical median consensus tree (d)}
    \label{fig:trees}
\end{figure}
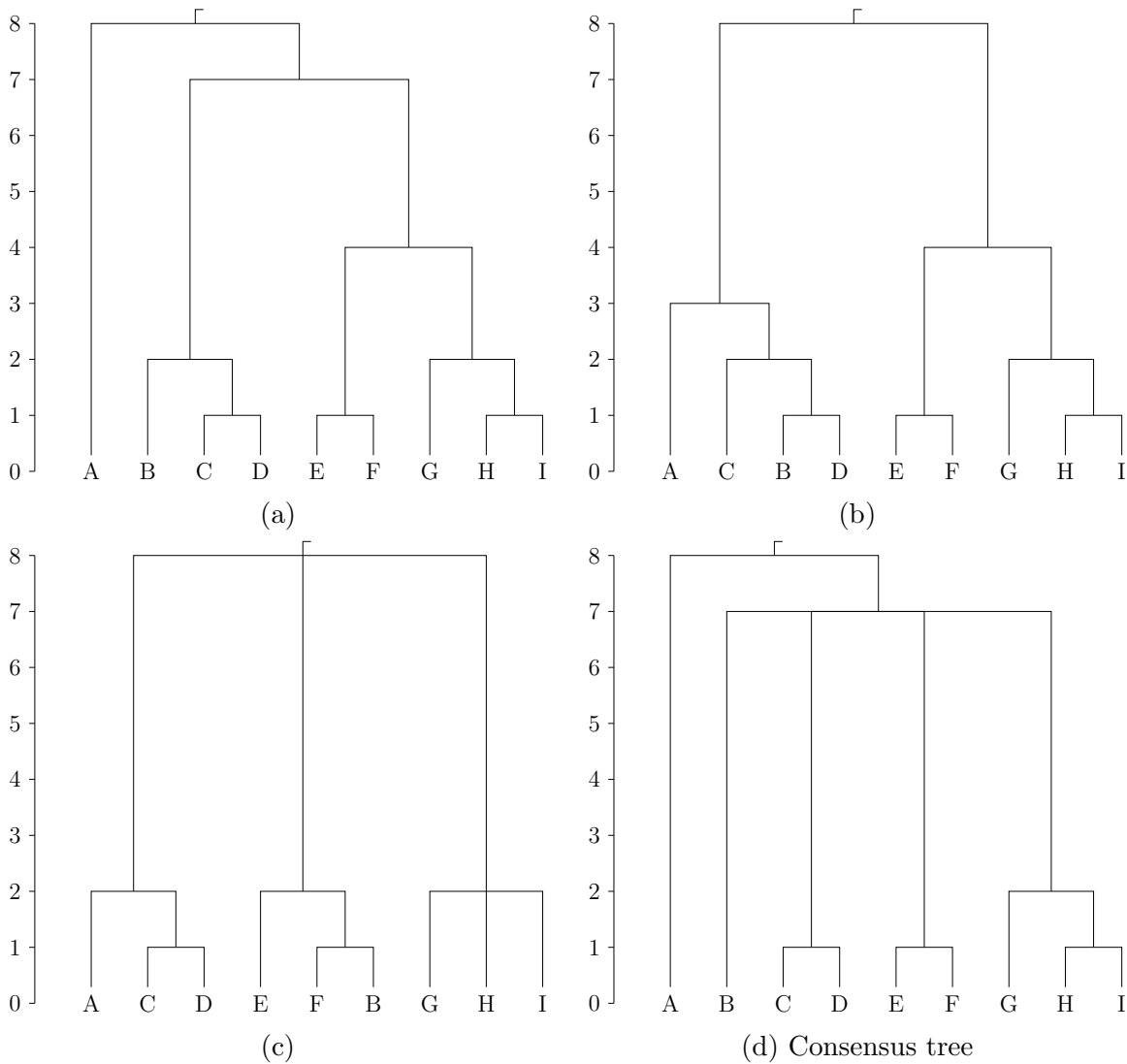

\subsection{Consensus trees}

Our goal now is to employ the results obtained so far to study the consensus problem for metric trees.
Formally, a \emph{consensus method} on equidistant trees, with $n$ taxa, is a function $c:\cT_n^*\to\cT_n$ where $\cT_n^*=\bigcup_{m\geq 1}\cT_n^m$.
For surveys on the subject see \cite{Bryant:2003} and, for a more recent account, \cite{Bryant+Francis+Steel:2017}.
We say that a consensus method $c$ is \emph{tropically convex} if $c(D_1,\dots,D_m)\in\tconv^{\max}(D_1,\dots,D_m)$ for every $m\geq 1$ and every $D_1,\dots,D_m\in\cT_n$.
Theorem \ref{thm:trop_consensus} says that selecting an arbitrary tree in $\FW(D_1,\dots,D_m)$ yields a consensus method which is tropically convex.
Recall that $\FW(D_1,\dots,D_m)$ is a polytrope in $\torus{\tbinom{n}{2}}$, which thus has at most $\tbinom{n}{2}$ tropical vertices.

\begin{definition}
  We define the \emph{tropical median consensus tree} of $D_1,\dots,D_m\in\cT_n$ as the ordinary average of the tropical vertices of $\FW(D_1,\dots,D_m)$.
\end{definition}

That ordinary average is the barycenter of the ordinary simplex spanned by the tropical vertices.
As polytropes are convex in the ordinary sense, the tropical median consensus tree method is tropically convex.
By Corollary~\ref{cor:strongPoly}, the tropical vertices of the Fermat--Weber set can be found in strongly polynomial time and hence also their ordinary average.

\begin{example}\label{exmp:newick}
  The first three equidistant trees on $n=9$ taxa in Fig.~\ref{fig:trees}, called (a), (b), (c), are taken from \cite[Chapter~7]{CladisticsBook:1998}.
  Since the trees in that reference are not metric, here we choose weights that are compatible with the graphical representation in \cite[Fig.~7.1]{CladisticsBook:1998}.
  The tropical median consensus tree is depicted in Fig.~\ref{fig:trees} as (d).
  This is the unique Fermat--Weber point of the three input trees in $\cT_9$.
  This can be verified using version~4.9 of \polymake \cite{DMV:polymake}.

  The \emph{Newick format} is a standard for encoding phylogenetic trees \cite[p.~590]{Felsenstein:2003}; it is supported by \polymake.
  A tree is represented as a string, where leaves are given by their (text) labels, and internal nodes correspond to matching pairs of parentheses.
  Recursively, such a pair of parentheses encloses the Newick representation of the subtree rooted at that internal node.
  Further, each node is followed by a numerical value, after a colon, and this denotes the length of the edge between the node and its parent.
  For example, the Newick representation of (a) is  \texttt{(A:8,((B:2,(C:1,D:1):1):5,((E:1,F:1):3,(G:2,(H:1,I:1):1):2):3):1)}.
\end{example}

Bryant, Francis and Steel \cite{Bryant+Francis+Steel:2017} impose three conditions for a consensus method to be \emph{regular}:
\begin{itemize}
\item[(U)] The consensus of any number of copies of the same tree, $T$, is $T$;
\item[(A)] the consensus does not depend on the ordering of the trees;
\item[(N)] permuting the taxa in the input trees results in the same permutation of the taxa in the consensus.
\end{itemize}
These properties are called unanimity, anonymity, and neutrality, respectively.
For the tropical median consensus all of them are immediate:
Unanimity is due to fact that $d_\triangle$ is definite; anonymity follows from the commutativity of the addition; neutrality is a consequence of the invariance of $d_\triangle$ under the action of the symmetric group.
It then follows from \cite[Theorem~3]{Bryant+Francis+Steel:2017} that the tropical median consensus is not \enquote{extension stable}.

Our next step is to investigate properties of arbitrary tropically convex consensus methods.
To this end,  let $i,j,k\in[n]$ be pairwise distinct taxa in some equidistant tree such that the lowest common ancestor of $i$ and $j$ is a proper descendant of the lowest common ancestor of $i$, $j$, and $k$.
Then we say that these taxa form a \emph{rooted triplet}, and we denote it by $ij|k$.
If $D=(d_{ij})$ is its ultrametric distance, then $ij|k$ is a rooted triplet if and only if $d_{ij}<d_{ik}$.
Note that the ultrametric property implies $d_{jk}=d_{ik}$, so we also have $d_{ij}<d_{jk}$.
We denote by $r(D)$ the set of rooted triplets of the tree.
A consensus method is called \emph{Pareto on rooted triplets} if $\bigcap_{\ell\in[m]}r(D_\ell)\subseteq r(D)$; it is called \emph{co-Pareto on rooted triplets} if $r(D)\subseteq\bigcup_{\ell\in[m]}r(D_\ell)$;
here $D_1,\dots,D_m$ correspond to the input trees, and $D$ represents the consensus tree.
These are desirable properties for consensus methods; see \cite[\textsection 3]{Bryant:2003}.

\begin{proposition} \label{prop:Pareto}
  Any tropically convex consensus method is Pareto and co-Pareto on rooted triplets.
\end{proposition}

\begin{proof}
  We consider the set of equidistant trees containing the rooted triplet $ij|k$, which we denote
  \[
    \cT_n(ij|k) \ := \ \SetOf{D\in\cT_n}{d_{ij}<d_{ik}} \enspace.
  \]
  The key observation is that this set is tropically convex: it arises as the intersection of $\cT_n$ with an open tropical halfspace.
  Note that, therefore, the complement in $\cT_n$ is also tropically convex.

  Now let $D_1,\dots,D_m$ be ultrametrics and $D$ any point in their max-tropical convex hull.
  We need to verify the Pareto and co-Pareto properties.
  If the rooted triplet $ij|k$ belongs to $\bigcap_{\ell}r(D_\ell)$, then $D_\ell\in\cT_n(ij|k)$ for every $\ell\in[m]$.
  As $\cT_n(ij|k)$ is tropically convex, we have  $D\in\cT_n(ij|k)$.
  Thus, $ij|k$ also belongs to $r(D)$, showing that a tropically convex consensus method is Pareto on rooted triplets.

  Conversely, if $ij|k$ does not belong to $\bigcup_{\ell\in[m]}r(D_i)$, then the input ultrametrics $D_1,\dots,D_m$ lie in the complement $\cT_n\setminus\cT_n(ij|k)$.
  Again, the latter set of tropically convex, and thus $D\notin\cT_n(ij|k)$.
  We infer that $ij|k\notin r(D)$, and we conclude that a tropically convex consensus method is co-Pareto on rooted triplets.
\end{proof}

Before we continue to study our tropical median consensus method, we look at other ideas.

\begin{figure}
  \centering
  \resizebox{7.5cm}{!}{\begin{tikzpicture}[sloped]
    \node (A) at (0,0) {A};
    \node (B) at (1,0) {B};
    \node (C) at (2,0) {C};
    \node (D) at (3,0) {D};
    \node (E) at (4,0) {E};
    \node (F) at (5,0) {F};
    \node (G) at (6,0) {G};
    \node (H) at (7,0) {H};
    \node (I) at (8,0) {I};

    \node (CD) at (2.5,2) {};
    \node (EF) at (4.5,2) {};
    \node (GHI) at (7,2) {};
    \node (all) at (3.5,8) {};
    \node (root) at (3.5,8.25) {};
    
    \draw  (A) |- (all.center);
    \draw  (B) |- (all.center);
    \draw  (C) |- (CD.center);
    \draw  (D) |- (CD.center);
    \draw  (E) |- (EF.center);
    \draw  (F) |- (EF.center);
    \draw  (G) |- (GHI.center);
    \draw  (H) |- (GHI.center);
    \draw  (I) |- (GHI.center);
    \draw  (CD.center) |- (all.center);
    \draw  (EF.center) |- (all.center);
    \draw  (GHI.center) |- (all.center);
    \draw  (all.center) |- (root);
    
    \draw (-1,0) -- (-1,8);
    
    \foreach \y in {0,1,2,3,4,5,6,7,8}                     
    \draw[shift={(0,\y)},color=black] (-1,0) -- (-1.1,0);
    
    \foreach \y in {0,1,2,3,4,5,6,7,8}   
    \node[left] at (-1.1,\y) {$\y$} ;
%
%
%
%
%
    \end{tikzpicture}}
  \caption{Pointwise maximum of the three trees from Example~\ref{exmp:trop_center}}
  \label{fig:trop_center}
\end{figure}
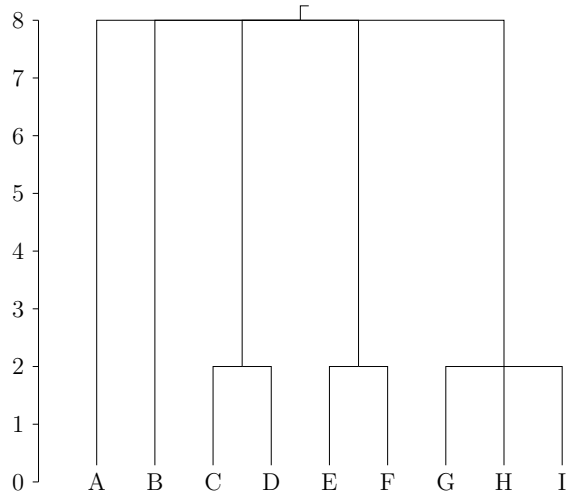
  
\begin{example} \label{exmp:trop_center}
  A particularly simple way to produce a tropically convex consensus tree method is the following.
  For given ultrametrics $D_1,\dots, D_m\in\cT_n$, we can consider the pointwise maximum $c(D_1,\dots,D_m)=D_1\oplus\dots\oplus D_m$.
  See Fig.~\ref{fig:trop_center} for the pointwise maximum of the three trees from Example~\ref{exmp:trop_center} and Fig.~\ref{fig:trees}.
  The tropical median consensus tree (d) in Fig.~\ref{fig:trees} partially resolves the maximum consensus tree from Fig.~\ref{fig:trop_center}.
	
  Note that the above definition depends on the representatives of $D_1,\dots,D_m$ modulo $\RR\1$.
  For the output in Fig.~\ref{fig:trop_center}, we used the representatives displayed in Fig.~\ref{fig:trees}.
  However, we could have chosen the representatives lying on the hyperplane $\cH$ in $\torus{\tbinom{n}{2}}$; the corresponding pointwise maximum represents the center of the smallest ball with respect to $d_\triangle$ that contains the points $D_1,\dots,D_m$.
  Alternatively, considering representatives with a fixed entry equal to 0, we obtain the \emph{tropical barycenter} from \cite[\S3.2]{ABGJ:2018}.

  Whatever convention we may fix for the representatives, the pointwise maximum exhibits a drawback.
  To exemplify it, consider the trees (a), (b), and one million copies of the tree (c) from Fig.~\ref{fig:trees}.
  Then, the pointwise maximum consensus is still the one from Fig.~\ref{fig:trop_center}, whereas the tropical median consensus tree looks like (c).
  So the pointwise maximum consensus is highly sensitive to outliers.
  In contrast, the tropical median consensus is robust.
\end{example}


Most known consensus tree methods deal with unweighted phylogenetic trees, so they may be seen as discrete analogues of our approach.
For instance, unweighted \enquote{median consensus trees} were defined by Bath\'el\'emy and McMorris \cite{MedianTree}, and the asymmetric case was analyzed by Phillips and Warnow \cite{AsymmMedianTree}.
Bryant \cite{Bryant:2003} presents only one consensus method for rooted trees that takes the edge lengths into consideration: the \enquote{average consensus tree} of Lapointe and Cucumel \cite{AvgConsensus}.
In \cite[\S2.4.2]{Bryant:2003} two drawbacks of the average consensus tree method are explicitly mentioned: no efficient algorithm is known, and the (co-)Pareto properties are unclear.
The average consensus method involves the Euclidean distance, and the unconstrained optima might lie outside the tree space $\cT_n$.
Therefore, ultrametric conditions must be imposed to the solution, making it difficult to obtain a regular consensus method; see \cite{AvgConsensus} for details.
Further, Lapointe and Cucumel \cite{Lapointe+Cucumel:2002} show that the procedure proposed is $\mathsf{NP}$-hard, and the solution may not be unique.
A similar complexity result exists for the median consensus method developed by Lavasseur and Lapointe \cite{Levasseur+Lapointe:2002}; see \cite{Day:87}.

In the remainder of this section we compare the tropical median consensus method to algorithms proposed by Lin and Yoshida \cite{Lin-Yoshida:2018}.
In fact, our approach is very similar and was inspired by that article.
Lin and co-authors \cite{ConvTreeSp} studied the Fermat--Weber problem for the symmetric tropical distance function, with a focus on the tree space.
Crucially, the symmetric tropical Fermat--Weber points may lie outside the $\max$-tropical convex hull; see \cite[Example~27]{ConvTreeSp}.

\begin{table}
  \caption{Timings (in seconds) for computing symmetric tropical Fermat--Weber sets of equidistant trees. The computations up to 8 leaves are averages over 10 iterations, whereas those for 9 and 10 leaves are results after one iteration}
  \label{table:sym-times}
  \centering
  \scalebox{0.75}{\begin{tabular}{crrrrrrrrrr}
  \toprule
    Leaves$\backslash$Trees & \multicolumn{1}{c}{1} & \multicolumn{1}{c}{2} & \multicolumn{1}{c}{3} & \multicolumn{1}{c}{4} & \multicolumn{1}{c}{5} & \multicolumn{1}{c}{6} & \multicolumn{1}{c}{7} & \multicolumn{1}{c}{8} & \multicolumn{1}{c}{9} & \multicolumn{1}{c}{10}  \\
    \midrule
    4 & 0.24 & 0.24 & 0.25 & 0.26 & 0.28 & 0.31 & 0.33 & 0.36 & 0.40 & 0.45 \\ 
    5 & 0.75 & 0.80 & 0.91 & 1.07 & 1.28 & 1.55 & 1.88 & 2.31 & 2.82 & 3.46 \\ 
    6 & 2.05 & 2.77 & 3.78 & 5.25 & 7.09 & 9.65 & 12.73 & 16.67 & 21.79 & 27.86 \\ 
    7 & 6.06 & 10.40 & 16.82 & 25.35 & 37.90 & 55.21 & 76.60 & 105.43 & 140.58 & 183.15 \\ 
    8 & 19.97 & 40.59 & 74.28 & 121.04 & 186.41 & 277.18 & 398.74 & 554.41 & 748.75 & 980.57 \\
    9 & 56.13 & 124.75 & 268.57 & 458.23 &  713.34 & 1094.16 & 1615.68 & 2291.03 &	3123.42 &  4186.37 \\
    10& 157.39 & 444.24 & 973.56 & 1705.29 & 2777.11 & 4127.49 & 6247.95 & 8874.29 & 12079.40 & 16134.13 \\
    \bottomrule
 \end{tabular}}
\end{table}

\begin{example}
  Symmetric tropical Fermat--Weber sets may be surprisingly complicated.
  Consider the trees
  \begin{equation}\label{eq:two-trees}
    T_1 = \mathtt{(D:10,(C:4,(B:2,A:2):2):6)} \quad\text{and}\quad T_2 = \mathtt{(A:10,(B:4,(C:2,D:2):2):6)}
  \end{equation}
  from \cite[Fig.~5]{Lin+Monod+Yoshida:2022}; here and below we employ the Newick format discussed in Example~\ref{exmp:newick}.
  The symmetric Fermat--Weber set of $T_1$ and $T_2$, denoted $\FWsym(T_1,T_2)$, contains the tropical segment between them, which exhibits seven distinct tree topologies, four of which are binary.
  In contrast, the asymmetric setting is trivial: the tropical median is $\mathtt{(A:10,D:10,(B:4,C:4):6)}$, and this is the unique asymmetric tropical Fermat--Weber point.
\end{example}
In view of \cite[Lemma~3.5]{Page+Yoshida+Zhang:2020} the symmetric tropical distance function might lead to a robust tropical consensus method.
However, examples like the above form a challenge to defining a method which is globally consistent.
The next case shows more differences between the symmetric and the asymmetric distances.

\begin{example}
  With $T_1$ and $T_2$ defined as in \eqref{eq:two-trees}, we consider two copies of $T_1$, two copies of $T_2$, and the tree $T_3=\mathtt{(A:10,((B:4,C:4):3,D:7):3)}$.
  That is to say, with multiplicities, we have five trees altogether.
  The unique asymmetric tropical Fermat--Weber point of these five trees is the tree $\mathtt{(A:10,D:10,(B:4,C:4):6)}$.
  On the other hand, $T_3$ is the unique symmetric tropical Fermat--Weber point, by \cite[Lemma~8]{Lin-Yoshida:2018}.
  Both Fermat--Weber points are unique, but they differ.
\end{example}

\begin{remark}
  The $\max$-tropical convexity of the tropical median consensus method ultimately rests on the specific formulation of the Fermat--Weber problem in \eqref{eq:from-the-sites}.
  Exchanging the arguments in the asymmetric tropical distance function gives the $\min$-tropical analog.
\end{remark}

\subsection{Computational experiments}

We compare running times for experiments concerning Fermat--Weber sets in tree space with respect to the symmetric and the asymmetric tropical distance functions.
As input data we take random trees which were produced using the function \texttt{rmtree} from the \texttt{R} library \texttt{ape} \cite{ape-software}.
This is similar to the experiment reported in \cite[Example 8]{YoshidaZhangZhang:2019}.
Most of the trees generated are not equidistant, so we adjust the lengths of the leaf edges to make them equidistant.
In this way we get any number of trees in $\cT_n$, for various values of $n$, the number of taxa.
Recall that, by \cite[Proposition 26]{ConvTreeSp}, the symmetric tropical Fermat--Weber set $\FWsym(T_1,\dots,T_m)$ of trees $T_i\in\cT_n$ is a convex polytope in $\torus{\tbinom{n}{2}}$.
The asymmetric tropical Fermat--Weber set $\FW(T_1,\dots,T_m)$ is a polytrope, and thus a polytope, too; cf.\ Theorem~\ref{thm:convFW}.
All timings are obtained with \polymake, version~4.9, running on a quad core Intel Core~i5-4590 processor (6599.89 bogomips), openSUSE Leap~15.3 (Linux~6.1.0).
For details see our data repository at \url{https://github.com/micjoswig/TropicalDataAnalysis/tree/main/Tropical_medians_by_transportation}.

\begin{table}
  \caption{Timings (in seconds) for computing asymmetric tropical Fermat--Weber sets of equidistant trees.
    Each entry is the average running time from 100 individual experiments, which show only very small variance}
\label{table:asym-times}
\centering
\begin{tabular}{crrrrrrrrrr}
 \toprule
 Leaves$\backslash$Trees & \multicolumn{1}{c}{1} & \multicolumn{1}{c}{2} & \multicolumn{1}{c}{3} & \multicolumn{1}{c}{4} & \multicolumn{1}{c}{5} & \multicolumn{1}{c}{6} & \multicolumn{1}{c}{7} & \multicolumn{1}{c}{8} & \multicolumn{1}{c}{9} & \multicolumn{1}{c}{10}  \\
 \midrule
 4 & 0.01 & 0.01 & 0.01 & 0.01 & 0.02 & 0.02 & 0.03 & 0.03 & 0.04 & 0.05 \\ 
 5 & 0.01 & 0.01 & 0.02 & 0.03 & 0.04 & 0.05 & 0.07 & 0.10 & 0.14 & 0.19 \\ 
 6 & 0.02 & 0.02 & 0.03 & 0.05 & 0.09 & 0.14 & 0.20 & 0.29 & 0.41 & 0.55 \\ 
 7 & 0.03 & 0.04 & 0.07 & 0.12 & 0.21 & 0.33 & 0.52 & 0.76 & 1.08 & 1.48 \\ 
 8 & 0.04 & 0.06 & 0.13 & 0.25 & 0.44 & 0.73 & 1.14 & 1.70 & 2.53 & 3.48 \\ 
 9 & 0.06 & 0.11 & 0.24 & 0.49 & 0.91 & 1.55 & 2.52 & 3.76 & 5.90 & 8.83 \\ 
 10 & 0.09 & 0.19 & 0.44 & 0.94 & 1.82 & 3.20 & 5.32 & 8.64 & 13.78 & 19.93 \\ 
\bottomrule 
\end{tabular} 
\end{table}

\paragraph{Entire Fermat--Weber sets.}
First, we compute exact facet descriptions of the polytopes $\FWsym(T_1,\dots,T_m)$ and $\FW(T_1,\dots,T_m)$, where $T_i\in\cT_n$, for various values of $m$ and $n$.
While the trees are generated with edge lengths given by floating point numbers, we convert them to exact rationals.
We are not aware of a published implementation for computing symmetric tropical Fermat--Weber sets, so we implemented it in \polymake.
The algorithm suggested by \cite[Proposition 26]{ConvTreeSp} allows for an improvement by exploiting properties of polyhedral L-convex functions in the sense of Murota \cite[\textsection 7.8]{Murota:2003}.
In this way, a facet description of $\FWsym(T_1,\dots,T_m)$ can be obtained from solving $n(n-1)$ linear programs similar to those in \eqref{eq:facetFW}.
The timings for up to ten trees on up to ten leaves are given in Table~\ref{table:sym-times}.
Lin and Yoshida report about computations of symmetric tropical Fermat--Weber sets in \cite[\S4]{Lin-Yoshida:2018}; however, they do not give timings.
The parameters leading to \cite[Table~1]{Lin-Yoshida:2018} are much smaller than the parameters in our Table~\ref{table:sym-times}.
The combinatorial description of the asymmetric tropical Fermat--Weber set $\FW(T_1,\dots,T_m)$ via optimal dual variables of a transportation problem allows us to exploit complementary slackness.
So we can compute the description of $\FW(T_1,\dots,T_m)$ in terms of its ordinary facets much faster than its symmetric counterpart; see Table~\ref{table:asym-times}.

\paragraph{Tropical median consensus trees.}
Our most interesting experiment is concerned with computing tropical median consensus trees of $m$ trees in $\cT_n$, again for various values of $m$ and $n$.
This time we use \mcf \cite{mcf} (through our \polymake interface), which is a standard implementation of the network simplex algorithm, using floating point computations.
Table~\ref{table:mcf-times} has the timings for up to 25 trees on up to 300 leaves.

The last row of Table~\ref{table:mcf-times}, for $m=25$ trees, is particularly interesting.
Observe that the timings in that row, with an increasing number of leaves, are not monotone.
The most likely explanation comes from Corollary~\ref{cor:trop_consensus}, which gives an upper bound for the dimension of $\FW(T_1,\dots,T_m)$; the critical contribution is $\gcd(m,\tbinom{n}{2})$.
We have $\tbinom{25}{2}=300$, and the running times are almost proportional to $\gcd(m,300)$.

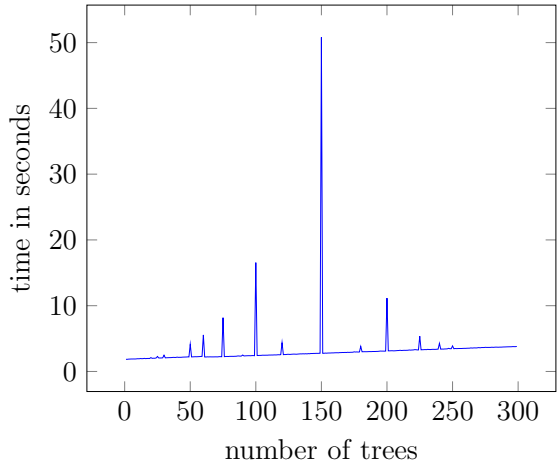
\begin{figure}[th]
    \centering
    \begin{tikzpicture}[scale=0.9]
			\begin{axis}[
			title=Tropical median consensus on trees with 25 leaves,
			xlabel={number of trees},
			ylabel={time in seconds},
			]
			\addplot [blue] table {25leaves.dat};
			\end{axis}
			\end{tikzpicture}
    \caption{Number of trees vs. time in seconds for up to 299 trees on 25 leaves}
    \label{fig:times_graph}
\end{figure}

In order to see what is going on, we conducted a more refined experiment for $n=25$ taxa, trying all values of $m$ in the set $\{1,2,\dots,299\}$.
The results are shown in Fig.~\ref{fig:times_graph}.
The timing for $m=300$, being more than 380 seconds, has been omitted for better readability of the other results.
There are only ten computations which last more than four seconds.
These are the values $m\in\{50,60,75,100,120,150,200,225,240,300\}$; all are integers with $\gcd(m,300)\geq 50$.
A more detailed analysis is beyond the scope of the present article.

To the best of our knowledge, no regular consensus method arising from the symmetric tropical distance function on tree space has been proposed. 
So there is no direct way to make a comparison.




\begin{table}[bh]
  \caption{Timings (in seconds) for computing tropical median consensus trees using \mcf \cite{mcf}}
  \label{table:mcf-times}
  \centering
  \begin{tabular}{crrrrrr}
  	\toprule
  	Leaves$\backslash$Trees  & \multicolumn{1}{c}{50} & \multicolumn{1}{c}{100} & \multicolumn{1}{c}{150} & \multicolumn{1}{c}{200} & \multicolumn{1}{c}{250} & \multicolumn{1}{c}{300}\\
  	\midrule
  	5 & 0.04 & 0.06 & 0.07 & 0.08 & 0.09 & 0.11 \\ 
  	10 & 0.11 & 0.16 & 0.23 & 0.26 & 0.31 & 0.36 \\ 
  	15 & 0.33 & 0.45 & 0.57 & 0.69 & 0.79 & 0.91 \\ 
  	20 & 0.87 & 1.08 & 1.29 & 1.50 & 1.70 & 1.92 \\
  	25 & 4.13 & 16.55 & 50.81 & 11.15 & 3.89 & 382.89 \\
  	\bottomrule 
  \end{tabular} 
\end{table}

\paragraph{Apicomplexa gene trees.}
For our final experiment we consider an existing dataset of $m=268$ trees with $n=8$ leaves, which was already studied by Page, Yoshida and Zhang \cite[\S6.2]{Page+Yoshida+Zhang:2020}.
Via simulated annealing, the latter authors generate three trees whose tropical convex hull fits best the input data; then they project the input trees on this tropical triangle and display the result in \cite[Fig.~6]{Page+Yoshida+Zhang:2020}.
So this line of research is more about dimension reduction techniques for analyzing data rather than obtaining location statistics; see also \cite{YoshidaZhangZhang:2019}.
In particular, this does not seem to lead to a regular consensus method in the sense of \cite{Bryant+Francis+Steel:2017}.
Nonetheless, our method applies to their data.

The trees discussed in \cite[\S6.2]{Page+Yoshida+Zhang:2020} have been reconstructed by Kuo, Wares and Kissinger \cite{apicomplexa} from 268 orthologous sequences with eight species of protozoa.
Seven species among the taxa are \textit{Babesia bovis} (Bb), \textit{Cryptosporidium parvum} (Cp), \textit{Eimeria tenella} (Et), \textit{Plasmodium falciparum} (Pf), \textit{Plasmodium vivax} (Pv), \textit{Theileria annulata} (Ta) and \textit{Toxoplasma gondii} (Tg).
The eighth taxon is \textit{Tetrahymena thermophila} (Tt), which forms the outgroup.
For this data we obtain the tropical median consensus tree
\[
  \begin{split}
    \texttt{(Cp:0.570333,}&\texttt{Et:0.570333,Tg:0.570333,Tt:0.570333,}\\
    &\texttt{(Pf:0.43862,Pv:0.43862):0.131713,(Bb:0.57033,Ta:0.57033):0.000003)} \enspace ,
  \end{split}
\]
which is displayed in Fig.~\ref{fig:api_tree}.
Only two cherries are resolved; the outgroup was not detected.
So our method is quite conservative, making an effort to avoid false positive results.
This may be an advantage, in particular since the 268 input trees were generated by a diverse range of methods.


\section{Conclusion}
\noindent
There is a considerable amount of work on the tropical metric geometry of the space of equidistant trees \cite{ConvTreeSp,Lin-Yoshida:2018,YoshidaZhangZhang:2019,Page+Yoshida+Zhang:2020,Lin+Monod+Yoshida:2022}.
This is quite natural, because the symmetric tropical distance between two equidistant trees can be interpreted naturally, as the cost of changing one tree into the other along a tropical line segment, where the cost is measured in the $\ell^\infty$ norm; see \cite{Lin+Monod+Yoshida:2022} for details.
Here we study the asymmetric tropical distance, which has a similar interpretation, but for the $\ell^1$ norm.
What makes the asymmetric tropical distance attractive is the fact that it leads to a regular consensus method, while this is not known to exist in the symmetric case.
Further, the tropical median consensus method is robust and fast in practice; it can be applied to hundreds of trees with dozens of taxa.
The robustness is computationally valuable as it makes floating-point computations reliable.

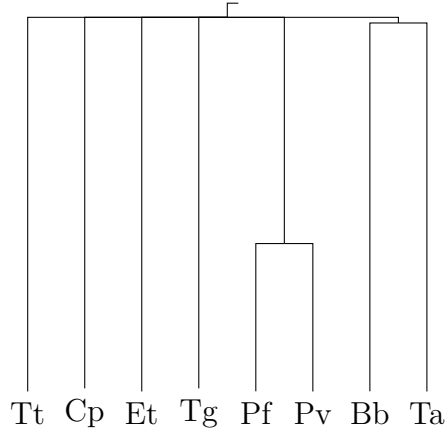
\begin{figure}
  \centering
  \begin{tikzpicture}[scale=0.75,sloped]
    \node (Tt) at (0,0) {Tt};
    \node (Cp) at (1,0) {Cp};
    \node (Et) at (2,0) {Et};
    \node (Tg) at (3,0) {Tg};
    \node (Pf) at (4,0) {Pf};
    \node (Pv) at (5,0) {Pv};
    \node (Bb) at (6,0) {Bb};
    \node (Ta) at (7,0) {Ta};
    
    \node (P) at (4.5,3) {};
    \node (BT) at (6.5,6.9) {};
    \node (all) at (3.5, 7) {};
    \node (root) at (3.5,7.25) {};
    
    \draw  (Tt) |- (all.center);
    \draw  (Cp) |- (all.center);
    \draw  (Et) |- (all.center);
    \draw  (Tg) |- (all.center);
    \draw  (Pf) |- (P.center);
    \draw  (Pv) |- (P.center);
    \draw  (Bb) |- (BT.center);
    \draw  (Ta) |- (BT.center);
    \draw  (P.center) |- (all.center);
    \draw  (BT.center) |- (all.center);
    \draw  (all.center) |- (root);
  \end{tikzpicture}
  \caption{The tropical median consensus of the apicomplexa data}
  \label{fig:api_tree}
\end{figure}

The tropical median consensus method seems to be rather conservative, in the sense that our consensus trees sometimes show little resolution.
This is not necessarily bad.
For instance, the apicomplexa gene trees studied in Section~\ref{sec:trees} come from a diverse mix of tree building methods, which should probably make it difficult to find a fully resolved consensus; cf.\ Fig.~\ref{fig:api_tree}.

\section*{Acknowledgements}
  We are indebted to Charles Semple and two anonymous reviewers for their comments on a previous version of this article.

%
%


 \printbibliography

@book{ETC,
  author = {Joswig, Michael},
  title = {Essentials of tropical combinatorics},
  publisher = {American Mathematical Society},
  address = {Providence, RI},
  series = {Graduate Studies in Mathematics},
  volume = {219},
  year = {2021},
}

@article {RRlattice,
    AUTHOR = {Amini, Omid and Manjunath, Madhusudan},
     TITLE = {Riemann-{R}och for sub-lattices of the root lattice {$A_n$}},
   JOURNAL = {Electron. J. Combin.},
  FJOURNAL = {Electronic Journal of Combinatorics},
    VOLUME = {17},
      YEAR = {2010},
    NUMBER = {1},
     PAGES = {Research Paper 124, 50},
   MRCLASS = {06B99 (05C50 05E99 52C07)},
  MRNUMBER = {2729373},
MRREVIEWER = {Tatiana Smirnova-Nagnibeda},
       __URL =
              {http://www.combinatorics.org/Volume_17/Abstracts/v17i1r124.html},
}

@article {ConvTreeSp,
    AUTHOR = {Lin, Bo and Sturmfels, Bernd and Tang, Xiaoxian and Yoshida, Ruriko},
     TITLE = {Convexity in tree spaces},
   JOURNAL = {SIAM J. Discrete Math.},
  FJOURNAL = {SIAM Journal on Discrete Mathematics},
    VOLUME = 31,
      YEAR = 2017,
    NUMBER = 3,
     PAGES = {2015--2038},
      ISSN = {0895-4801},
   MRCLASS = {05C90 (05C05 05C12 14T05 30F45 52B40 68U05 92D15)},
  MRNUMBER = 3693600,
MRREVIEWER = {Mareike Fischer},
       DOI = {10.1137/16M1079841},
       _URL = {https://doi.org/10.1137/16M1079841},
}

@article {Lin-Yoshida:2018,
    AUTHOR = {Lin, Bo and Yoshida, Ruriko},
     TITLE = {Tropical {F}ermat-{W}eber points},
   JOURNAL = {SIAM J. Discrete Math.},
  FJOURNAL = {SIAM Journal on Discrete Mathematics},
    VOLUME = {32},
      YEAR = {2018},
    NUMBER = {2},
     PAGES = {1229--1245},
      ISSN = {0895-4801},
   MRCLASS = {13P25 (14T90 52B11 92B15)},
  MRNUMBER = {3810501},
MRREVIEWER = {Mateusz Micha\l ek},
       DOI = {10.1137/16M1071122},
       _URL = {https://doi.org/10.1137/16M1071122},
}

@article {ArdilaKlivans:2006,
    AUTHOR = {Ardila, Federico and Klivans, Caroline J.},
     TITLE = {The {B}ergman complex of a matroid and phylogenetic trees},
   JOURNAL = {J. Combin. Theory Ser. B},
  FJOURNAL = {Journal of Combinatorial Theory. Series B},
    VOLUME = {96},
      YEAR = {2006},
    NUMBER = {1},
     PAGES = {38--49},
      ISSN = {0095-8956},
   MRCLASS = {05B35},
  MRNUMBER = {2185977},
MRREVIEWER = {Neil L. White},
       DOI = {10.1016/j.jctb.2005.06.004},
       _URL = {https://doi.org/10.1016/j.jctb.2005.06.004},
}

@article {YoshidaZhangZhang:2019,
    AUTHOR = {Yoshida, Ruriko and Zhang, Leon and Zhang, Xu},
     TITLE = {Tropical principal component analysis and its application to phylogenetics},
   JOURNAL = {Bull. Math. Biol.},
  FJOURNAL = {Bulletin of Mathematical Biology. A Journal Devoted to
              Research at the Interface of the Life and Mathematical
              Sciences},
    VOLUME = {81},
      YEAR = {2019},
    NUMBER = {2},
     PAGES = {568--597},
      ISSN = {0092-8240},
   MRCLASS = {14T90 (92D15)},
  MRNUMBER = {3902912},
MRREVIEWER = {Boulos El Hilany},
       DOI = {10.1007/s11538-018-0493-4},
       _URL = {https://doi.org/10.1007/s11538-018-0493-4},
}

@article {LapLatticeGraph,
    AUTHOR = {Manjunath, Madhusudan},
     TITLE = {The {L}aplacian lattice of a graph under a simplicial distance function},
   JOURNAL = {European J. Combin.},
  FJOURNAL = {European Journal of Combinatorics},
    VOLUME = {34},
      YEAR = {2013},
    NUMBER = {6},
     PAGES = {1051--1070},
      ISSN = {0195-6698},
   MRCLASS = {52B20 (05C50 06A07)},
  MRNUMBER = {3037988},
MRREVIEWER = {Stephen J. Young},
       DOI = {10.1016/j.ejc.2013.01.010},
       __URL = {https://doi.org/10.1016/j.ejc.2013.01.010},
}

@book {BMS:1999,
    AUTHOR = {Boltyanski, Vladimir and Martini, Horst and Soltan, Valeriu},
     TITLE = {Geometric methods and optimization problems},
    SERIES = {Combinatorial Optimization},
    VOLUME = {4},
 PUBLISHER = {Kluwer Academic Publishers, Dordrecht},
      YEAR = {1999},
      ISBN = {0-7923-5454-0},
   MRCLASS = {90-02 (49K15 65D18 90B85)},
  MRNUMBER = {1677397},
MRREVIEWER = {Anita Kripfganz},
       DOI = {10.1007/978-1-4615-5319-9},
       _URL = {https://doi.org/10.1007/978-1-4615-5319-9},
}

@book {DLRS,
    AUTHOR = {De Loera, Jes\'{u}s A. and Rambau, J\"{o}rg and Santos, Francisco},
     TITLE = {Triangulations},
    SERIES = {Algorithms and Computation in Mathematics},
    VOLUME = {25},
      NOTE = {Structures for algorithms and applications},
 PUBLISHER = {Springer-Verlag, Berlin},
      YEAR = {2010},
      ISBN = {978-3-642-12970-4},
   MRCLASS = {52B55 (05C10 52B05 57Q15 68U05)},
  MRNUMBER = {2743368},
       DOI = {10.1007/978-3-642-12971-1},
       _URL = {https://doi.org/10.1007/978-3-642-12971-1},
}

@book {Schrijver:1986,
    AUTHOR = {Schrijver, Alexander},
     TITLE = {Theory of linear and integer programming},
    SERIES = {Wiley-Interscience Series in Discrete Mathematics},
      NOTE = {A Wiley-Interscience Publication},
 PUBLISHER = {John Wiley \& Sons, Ltd., Chichester},
      YEAR = {1986},
      ISBN = {0-471-90854-1},
   MRCLASS = {90C05 (90C10)},
  MRNUMBER = {874114},
MRREVIEWER = {J\"{u}rgen K\"{o}hler},
}

@article {Brenner:2008,
    AUTHOR = {Brenner, Ulrich},
     TITLE = {A faster polynomial algorithm for the unbalanced {H}itchcock transportation problem},
   JOURNAL = {Oper. Res. Lett.},
  FJOURNAL = {Operations Research Letters},
    VOLUME = {36},
      YEAR = {2008},
    NUMBER = {4},
     PAGES = {408--413},
      ISSN = {0167-6377},
   MRCLASS = {90B10 (68W40 90B06 90C35)},
  MRNUMBER = {2437262},
       DOI = {10.1016/j.orl.2008.01.011},
       _URL = {https://doi.org/10.1016/j.orl.2008.01.011},
}

@incollection {DeLoera+Kim:2014,
    AUTHOR = {De Loera, Jes\'{u}s A. and Kim, Edward D.},
     TITLE = {Combinatorics and geometry of transportation polytopes: an
              update},
 BOOKTITLE = {Discrete geometry and algebraic combinatorics},
    SERIES = {Contemp. Math.},
    VOLUME = {625},
     PAGES = {37--76},
 PUBLISHER = {Amer. Math. Soc., Providence, RI},
      YEAR = {2014},
   MRCLASS = {52B12 (52B05 90C08)},
  MRNUMBER = {3289405},
MRREVIEWER = {Alexander I. Barvinok},
       DOI = {10.1090/conm/625/12491},
       _URL = {https://doi.org/10.1090/conm/625/12491},
}

@article{Page+Yoshida+Zhang:2020,
    author = {Page, Robert and Yoshida, Ruriko and Zhang, Leon},
    title = {Tropical principal component analysis on the space of phylogenetic trees},
    journal = {Bioinformatics},
    volume = {36},
    number = {17},
    pages = {4590-4598},
    year = {2020},
    month = {06},
    issn = {1367-4803},
    doi = {10.1093/bioinformatics/btaa564},
    url = {https://doi.org/10.1093/bioinformatics/btaa564},
    __eprint = {https://academic.oup.com/bioinformatics/article-pdf/36/17/4590/34220689/btaa564.pdf},
}

@Unpublished{GaertnerJaggi:2006,
  author = 	 {G\"{a}rtner, Bernd and Jaggi, Martin},
  title = 	 {Tropical support vector machines},
  note = 	 {ACS Technical Report No.: ACS-TR-362502-01},
  year = 	 2006
}

@BOOK{ASCB,
  title = {Algebraic statistics for computational biology},
  publisher = {Cambridge University Press},
  year = {2005},
  editor = {Pachter, Lior and Sturmfels, Bernd},
  address = {New York},
  isbn = {978-0-521-85700-0; 0-521-85700-7},
  mrclass = {92-06 (62-06)},
  mrnumber = {2205865 (2006i:92002)}
}

@ARTICLE{BilleraHolmesVogtmann01,
  author = {Billera, Louis J. and Holmes, Susan P. and Vogtmann, Karen},
  title = {Geometry of the space of phylogenetic trees},
  journal = {Adv. in Appl. Math.},
  year = {2001},
  volume = {27},
  pages = {733--767},
  number = {4},
  doi = {10.1006/aama.2001.0759},
  fjournal = {Advances in Applied Mathematics},
  issn = {0196-8858},
  mrclass = {05C90 (92B10)},
  mrnumber = {1867931 (2002k:05229)},
  mrreviewer = {Charles A. Semple},
}

@article {Tokuyama+Nakano:1995,
    AUTHOR = {Tokuyama, Takeshi and Nakano, Jun},
     TITLE = {Efficient algorithms for the {H}itchcock transportation
              problem},
   JOURNAL = {SIAM J. Comput.},
  FJOURNAL = {SIAM Journal on Computing},
    VOLUME = {24},
      YEAR = {1995},
    NUMBER = {3},
     PAGES = {563--578},
      ISSN = {0097-5397},
   MRCLASS = {90B80 (68U05 90C08)},
  MRNUMBER = {1333855},
       DOI = {10.1137/S0097539792236729},
       _URL = {https://doi.org/10.1137/S0097539792236729},
}

@BOOK{Schrijver03:CO_A,
  title = {Combinatorial optimization. {P}olyhedra and efficiency. {V}ol. {A}},
  publisher = {Springer-Verlag},
  year = {2003},
  author = {Schrijver, Alexander},
  volume = {24},
  series = {Algorithms and Combinatorics},
  address = {Berlin},
  note = {Paths, flows, matchings, Chapters 1--38},
  isbn = {3-540-44389-4},
  mrclass = {90-02 (05-02 52B55 68Q25 68R10 90C27 90C35 90C57)},
  mrnumber = {MR1956924 (2004b:90004a)},
  mrreviewer = {Alexander I. Barvinok}
}

@article {Kleinschmidt+Schannath:1995,
    AUTHOR = {Kleinschmidt, Peter and Schannath, Heinz},
     TITLE = {A strongly polynomial algorithm for the transportation problem},
   JOURNAL = {Math. Programming},
  FJOURNAL = {Mathematical Programming},
    VOLUME = {68},
      YEAR = {1995},
    NUMBER = {1, Ser. A},
     PAGES = {1--13},
      ISSN = {0025-5610},
   MRCLASS = {90C08},
  MRNUMBER = {1312103},
MRREVIEWER = {Ulrich Faigle},
       DOI = {10.1007/BF01585755},
       _URL = {https://doi.org/10.1007/BF01585755},
}

@book {SempleSteel:2003,
    AUTHOR = {Semple, Charles and Steel, Mike},
     TITLE = {Phylogenetics},
    SERIES = {Oxford Lecture Series in Mathematics and its Applications},
    VOLUME = {24},
 PUBLISHER = {Oxford University Press, Oxford},
      YEAR = {2003},
      ISBN = {0-19-850942-1},
   MRCLASS = {92D15 (05C05 05C90 92D40)},
  MRNUMBER = {2060009},
MRREVIEWER = {Vincent L. Moulton},
}

@incollection {Bryant:2003,
    AUTHOR = {Bryant, David},
     TITLE = {A classification of consensus methods for phylogenetics},
 BOOKTITLE = {Bioconsensus ({P}iscataway, {NJ}, 2000/2001)},
    SERIES = {DIMACS Ser. Discrete Math. Theoret. Comput. Sci.},
    VOLUME = {61},
     PAGES = {163--183},
 PUBLISHER = {Amer. Math. Soc., Providence, RI},
      YEAR = {2003},
   MRCLASS = {92B10 (05C05 91B12 92D15)},
  MRNUMBER = {1982426},
}

@book {DantzigThapa:1997,
    AUTHOR = {Dantzig, George B. and Thapa, Mukund N.},
     TITLE = {Linear programming. 1},
    SERIES = {Springer Series in Operations Research},
 PUBLISHER = {Springer-Verlag, New York},
      YEAR = {1997},
      ISBN = {0-387-94833-3},
   MRCLASS = {90C05 (90-02)},
  MRNUMBER = {1485773},
MRREVIEWER = {K. G. Murty},
}

@article {BlockYu:2006,
    AUTHOR = {Block, Florian and Yu, Josephine},
     TITLE = {Tropical convexity via cellular resolutions},
   JOURNAL = {J. Algebraic Combin.},
  FJOURNAL = {Journal of Algebraic Combinatorics. An International Journal},
    VOLUME = {24},
      YEAR = {2006},
    NUMBER = {1},
     PAGES = {103--114},
      ISSN = {0925-9899},
   MRCLASS = {52B70 (13P99 14P99 52B10 52B55)},
  MRNUMBER = {2245783},
MRREVIEWER = {Domenico Fiorenza},
       DOI = {10.1007/s10801-006-9104-9},
       _URL = {https://doi.org/10.1007/s10801-006-9104-9},
}

@book {CladisticsBook:1998,
    AUTHOR = {Kitching, Ian J. and Forey, Peter L. and Humphries, Christopher J. and Williams, David M.},
     TITLE = {Cladistics: The Theory and Practice of Parsimony Analysis},
      EDITION = {2nd edition},
 PUBLISHER = {Oxford University Press},
      YEAR = {1998},
      ISBN = {019 850138 2},
}

@book {Tropical+Book,
    AUTHOR = {Maclagan, Diane and Sturmfels, Bernd},
     TITLE = {Introduction to tropical geometry},
    SERIES = {Graduate Studies in Mathematics},
    VOLUME = {161},
 PUBLISHER = {American Mathematical Society, Providence, RI},
      YEAR = {2015},
      ISBN = {978-0-8218-5198-2},
   MRCLASS = {14T05 (05B35 15A80 52B70)},
  MRNUMBER = {3287221},
}

@ARTICLE{Develin+Sturmfels:2004,
  author = {Develin, Mike and Sturmfels, Bernd},
  title = {Tropical convexity},
  journal = {Doc. Math.},
  year = {2004},
  volume = {9},
  pages = {1--27 (electronic)},
  note = {correction: ibid., pp.\ 205--206},
  fjournal = {Documenta Mathematica},
  issn = {1431-0635},
  mrclass = {52A30 (52B10 92D15)},
  mrnumber = {MR2054977 (2005i:52010)},
  mrreviewer = {Gerard Sierksma}
}

@book {Butkovic:2010,
    AUTHOR = {Butkovi{\v{c}}, Peter},
     TITLE = {Max-linear systems: theory and algorithms},
    SERIES = {Springer Monographs in Mathematics},
 PUBLISHER = {Springer-Verlag London, Ltd., London},
      YEAR = {2010},
      ISBN = {978-1-84996-298-8},
   MRCLASS = {15A80 (90C27 91B02 93B03)},
  MRNUMBER = {2681232 (2011e:15049)},
       DOI = {10.1007/978-1-84996-299-5},
}

@Article{MedianTree,
 Author = {Barth\'el\'emy, Jean-Pierre and McMorris, Frederick R.},
 Title = {The median procedure for $n$-trees},
 FJournal = {Journal of Classification},
 Journal = {J. Classif.},
 ISSN = {0176-4268},
 Volume = {3},
 Pages = {329--334},
 Year = {1986},
 Language = {English},
 DOI = {10.1007/BF01894194},
 zbMATH = {4001233},
 Zbl = {0617.62066}
}

@article {AsymmMedianTree,
    AUTHOR = {Phillips, Cynthia and Warnow, Tandy J.},
     TITLE = {The asymmetric median tree---a new model for building consensus trees},
   JOURNAL = {Discrete Appl. Math.},
  FJOURNAL = {Discrete Applied Mathematics. The Journal of Combinatorial Algorithms, Informatics and Computational Sciences},
    VOLUME = {71},
      YEAR = {1996},
    NUMBER = {1-3},
     PAGES = {311--335},
      ISSN = {0166-218X},
   MRCLASS = {92D25 (68Q25 92D15)},
  MRNUMBER = {1420306},
       DOI = {10.1016/S0166-218X(96)00071-6},
       _URL = {https://doi.org/10.1016/S0166-218X(96)00071-6},
}

@article{AvgConsensus,
 author = {Lapointe, Francois-Joseph and Cucumel, Guy},
 title = {The Average Consensus Procedure: Combination of Weighted Trees Containing Identical or Overlapping Sets of Taxa},
 ISSN = {10635157, 1076836X},
 doi = {10.2307/2413625},
 __URL = {http://www.jstor.org/stable/2413625},
 journal = {Systematic Biology},
 number = {2},
 pages = {306--312},
 publisher = {Oxford University Press, Society of Systematic Biologists},
 urldate = {2022-04-22},
 volume = {46},
 year = {1997}
}

@article{Bryant+Francis+Steel:2017,
    author = {Bryant, David and Francis, Andrew and Steel, Mike},
    title = "{Can We \enquote{Future-Proof} Consensus Trees?}",
    journal = {Systematic Biology},
    volume = {66},
    number = {4},
    pages = {611-619},
    year = {2017},
    month = {02},
    abstract = "{Consensus methods are widely used for combining phylogenetic trees into a single estimate of the evolutionary tree for a group of species. As more taxa are added, the new source trees may begin to tell a different evolutionary story when restricted to the original set of taxa. However, if the new trees, restricted to the original set of taxa, were to agree exactly with the earlier trees, then we might hope that their consensus would either agree with or resolve the original consensus tree. In this article, we ask under what conditions consensus methods exist that are “future proof” in this sense. While we show that some methods (e.g., Adams consensus) have this property for specific types of input, we also establish a rather surprising “no-go” theorem: there is no “reasonable” consensus method that satisfies the future-proofing property in general. We then investigate a second notion of “future proofing” for consensus methods, in which trees (rather than taxa) are added, and establish some positive and negative results. We end with some questions for future work.}",
    issn = {1063-5157},
    doi = {10.1093/sysbio/syx030},
    url = {https://doi.org/10.1093/sysbio/syx030},
    _eprint = {https://academic.oup.com/sysbio/article-pdf/66/4/611/17726420/syx030.pdf},
}

@article{ABGJ:2018,
  author = {Allamigeon, Xavier and Benchimol, Pascal and Gaubert, St\'ephane and Joswig, Michael},
  title =  {Log-barrier interior point methods are not strongly polynomial},
  journal = {SIAM J. Appl. Algebra Geom.},
  year = {2018},
  volume = 2,
  number = 1,
  pages = {140-178},
  doi = {10.1137/17M1142132},
  arxiv = {1708.01544}
}

@article {Lin+Monod+Yoshida:2022,
    AUTHOR = {Lin, Bo and Monod, Anthea and Yoshida, Ruriko},
     TITLE = {Tropical geometric variation of tree shapes},
   JOURNAL = {Discrete Comput. Geom.},
  FJOURNAL = {Discrete \& Computational Geometry. An International Journal
              of Mathematics and Computer Science},
    VOLUME = {68},
      YEAR = {2022},
    NUMBER = {3},
     PAGES = {817--849},
      ISSN = {0179-5376},
   MRCLASS = {14T90 (62R01 92-10)},
  MRNUMBER = {4481323},
       DOI = {10.1007/s00454-022-00410-y},
       _URL = {https://doi.org/10.1007/s00454-022-00410-y},
}

@incollection {Lapointe+Cucumel:2002,
    AUTHOR = {Lapointe, Fran\c{c}ois-Joseph and Cucumel, Guy},
     TITLE = {Multiple consensus trees},
 BOOKTITLE = {Classification, clustering, and data analysis ({C}racow, 2002)},
    SERIES = {Stud. Classification Data Anal. Knowledge Organ.},
     PAGES = {359--364},
 PUBLISHER = {Springer, Berlin},
      YEAR = {2002},
   MRCLASS = {62H30 (62P10)},
  MRNUMBER = {2010472},
}

@incollection {Levasseur+Lapointe:2002,
    AUTHOR = {Levasseur, Claudine and Lapointe, Fran\c{c}ois-Joseph},
     TITLE = {A family of average consensus methods for weighted trees},
 BOOKTITLE = {Classification, clustering, and data analysis ({C}racow,
              2002)},
    SERIES = {Stud. Classification Data Anal. Knowledge Organ.},
     PAGES = {365--369},
 PUBLISHER = {Springer, Berlin},
      YEAR = {2002},
   MRCLASS = {62H30},
  MRNUMBER = {2010473},
}

@article {Day:87,
    AUTHOR = {Day, William H. E.},
     TITLE = {Computational complexity of inferring phylogenies from
              dissimilarity matrices},
   JOURNAL = {Bull. Math. Biol.},
  FJOURNAL = {Bulletin of Mathematical Biology},
    VOLUME = {49},
      YEAR = {1987},
    NUMBER = {4},
     PAGES = {461--467},
      ISSN = {0092-8240},
   MRCLASS = {92A10 (68Q25 92A05)},
  MRNUMBER = {908160},
       DOI = {10.1016/S0092-8240(87)80007-1},
       _URL = {https://doi.org/10.1016/S0092-8240(87)80007-1},
}

@incollection {DMV:polymake,
    AUTHOR = {Gawrilow, Ewgenij and Joswig, Michael},
     TITLE = {\polymake: a framework for analyzing convex polytopes},
 BOOKTITLE = {Polytopes---combinatorics and computation (Oberwolfach, 1997)},
    SERIES = {DMV Sem.},
    VOLUME = {29},
     PAGES = {43--73},
 PUBLISHER = {Birk\-h\"au\-ser},
   ADDRESS = {Basel},
      YEAR = {2000},
   MRCLASS = {52B55 (68U05)},
  MRNUMBER = {MR1785292 (2001f:52033)},
}

@Article{ape-software,
    title = {ape 5.0: an environment for modern phylogenetics and
      evolutionary analyses in {R}},
    author = {E. Paradis and K. Schliep},
    journal = {Bioinformatics},
    year = {2019},
    volume = {35},
    pages = {526-528},
}

@book {Murota:2003,
    AUTHOR = {Murota, Kazuo},
     TITLE = {Discrete convex analysis},
    SERIES = {SIAM Monographs on Discrete Mathematics and Applications},
 PUBLISHER = {Society for Industrial and Applied Mathematics (SIAM), Philadelphia, PA},
      YEAR = {2003},
     PAGES = {xxii+389},
      ISBN = {0-89871-540-7},
   MRCLASS = {90-02 (52-02 90C27 90C46 91B02)},
  MRNUMBER = {1997998},
MRREVIEWER = {Ulrich Faigle},
       DOI = {10.1137/1.9780898718508},
       URL = {https://doi.org/10.1137/1.9780898718508},
}

@software{mcf,
  author = {L\"{o}bel, Andreas},
  title = {{MCF} -- A network simplex implementation},
  url = {https://www.zib.de/opt-long_projects/Software/Mcf/},
  version = {1.3},
  date = {2004},
}

@book{Felsenstein:2003,
  author = {Felsenstein, Joseph},
  publisher = {Sinauer Associates},
  title = {Inferring phylogenies},
  year = {2003}
}

@article{apicomplexa,
    author = {Kuo, Chih-Horng and Wares, John P. and Kissinger, Jessica C.},
    title = "{The Apicomplexan Whole-Genome Phylogeny: An Analysis of Incongruence among Gene Trees}",
    journal = {Molecular Biology and Evolution},
    volume = {25},
    number = {12},
    pages = {2689-2698},
    year = {2008},
    month = {09},
    issn = {0737-4038},
    doi = {10.1093/molbev/msn213},
    url = {https://doi.org/10.1093/molbev/msn213},
    eprint = {https://academic.oup.com/mbe/article-pdf/25/12/2689/13639732/msn213.pdf},
}

@Book{NickelPuerto:2005,
  author = 	 {Nickel, Stefan and Puerto Albandoz, Justo},
  title = 	 {Location Theory, A Unified Approach},
  publisher = 	 {Springer},
  year = 	 2005,
  doi = {https://doi.org/10.1007/3-540-27640-8},
}

\end{document}